\documentclass{gtpart}
\usepackage{amscd,amssymb,verbatim,graphicx,diagrams}
\usepackage[matrix,arrow,curve]{xy}
\usepackage{mathrsfs}
\usepackage{tikz}
\usetikzlibrary{arrows,decorations.markings}

\let\oldmarginpar\marginpar
\renewcommand\marginpar[1]{\-\oldmarginpar[\raggedleft\footnotesize #1]%
	{\raggedright\footnotesize #1}}

\newskip\stdskip                      % standard vertical space
\newcommand{\ssl}{{\mathfrak sl}}

\renewcommand{\mod}{\operatorname{mod}}

\newcommand{\OO}{\mathcal{O}}

\newcommand{\pr}{\operatorname{pr}}

\newcommand{\coker}{\operatorname{coker}}

\newcommand{\NN}{{\cal N}}

\newcommand{\be}{{\bf e}}

\renewcommand{\Re}{\operatorname{Re}}

\newcommand{\m}{\mathfrak{m}} 
\renewcommand{\k}{\mathbf{k}}

\newcommand{\lan}{\langle}
\newcommand{\ran}{\rangle}

\newcommand{\CC}{\mathscr{C}}
\newcommand{\UU}{\mathcal{U}}

\newcommand{\Mat}{\operatorname{Mat}}

\newcommand{\ga}{\gamma}
\newcommand{\de}{\delta}

\numberwithin{equation}{section}

\newcommand{\Perf}{\operatorname{Perf}}

\newcommand{\GL}{\operatorname{GL}}

\newcommand{\bm}{{\bf m}}
\newcommand{\T}{\mathbb{T}}

\newtheorem{thm}{Theorem}[subsection]
\newtheorem{prop}[thm]{Proposition}
\newtheorem{lem}[thm]{Lemma}
\newtheorem{cor}[thm]{Corollary}
{  \theoremstyle{definition}
\newtheorem{defi}[thm]{Definition}
\newtheorem{ex}[thm]{Example}

\newtheorem{rem}[thm]{Remark}

}

\newcommand{\Pf}{\noindent {\it Proof}}
\newcommand{\id}{\operatorname{id}}

\newcommand{\ov}{\overline}

\renewcommand{\Im}{\operatorname{Im}}

\renewcommand{\AA}{{\cal A}}

\newcommand{\FF}{{\cal F}}

\newcommand{\MM}{\mathcal{M}}
\newcommand{\TT}{\mathcal{T}}
\newcommand{\XX}{{\mathcal X}}
\newcommand{\YY}{{\mathcal Y}}

\newcommand{\ABD}{\operatorname{ABD}}

\newcommand{\Om}{\Omega}

\newcommand{\Hom}{\operatorname{Hom}}
\newcommand{\sHom}{{\mathcal H}om}

\newcommand{\Ext}{\operatorname{Ext}}

\renewcommand{\a}{\alpha}
\renewcommand{\b}{\beta}

\newcommand{\De}{\Delta}
\newcommand{\la}{\lambda}
\renewcommand{\th}{\theta}
\renewcommand{\C}{{\Bbb C}}

\renewcommand{\R}{{\Bbb R}}
\renewcommand{\Z}{{\Bbb Z}}

\newcommand{\wt}{\widetilde}
\newcommand{\ot}{\otimes}

\newcommand{\sub}{\subset}
\newcommand{\ed}{\qed\vspace{1mm}}

\renewcommand{\fg}{{\mathfrak g}}

\renewcommand{\P}{\mathbb{P}}
\newcommand{\rmP}{\mathrm{P}}

\title{Associative Yang-Baxter equation and\\ Fukaya categories of square-tiled surfaces}
\author[Yank{\i} Lekili and Alexander Polishchuk]{Yank{\i} Lekili \\ Alexander Polishchuk}

\address{King's College London}
\address{University of Oregon}

\begin{document}

\begin{abstract} We show that all strongly non-degenerate trigonometric solutions of the
    associative Yang-Baxter equation (AYBE) can be obtained from triple Massey products in the Fukaya
    categories of square-tiled surfaces. Along the way, we give a classification result for
    cyclic $A_\infty$-algebra structures on a certain Frobenius algebra associated with a pair of
    $1$-spherical objects in terms of the equivalence classes of the corresponding solutions of the
    AYBE. As an application, combining our results with homological
    mirror symmetry for punctured tori (cf. \cite{LP}),  we prove that any two simple vector bundles on a cycle of projective
    lines are related by a sequence of $1$-spherical twists and their inverses.
\end{abstract}

\maketitle

\section*{Introduction}

The associative Yang-Baxter equation (AYBE) is the equation
\begin{equation}\label{AYBE}
r^{12}(-u',v)r^{13}(u+u',v+v')-
r^{23}(u+u',v')r^{12}(u,v)+
r^{13}(u,v+v')r^{23}(u',v')=0,
\end{equation}
where $r: \C^2 \to \Mat_n(\C) \otimes \Mat_n(\C)$ is a meromorphic function of two complex variables $(u,v)$
in a neighborhood of $(0,0)$ taking values in $\Mat_n(\C) \otimes \Mat_n(\C)$, and $\Mat_n(\C)$ is
the matrix algebra. Here we use the notation $r^{12}=r\ot 1\in \Mat_n(\C) \ot \Mat_n(\C) \ot \Mat_n(\C)$, etc.
The equation \eqref{AYBE} is usually coupled with the \emph{skew-symmetry} (also called \emph{unitarity}) condition
\begin{equation}\label{AYBE-skew-eq}
r^{21}(-u,-v)=-r(u,v)
\end{equation}
where $r^{21}$ is obtained from $r$ by the transposition of tensor factors $a_2 \otimes a_1 \to a_1 \otimes a_2$. 
Note that the constant version of the AYBE was studied by Aguiar \cite{Aguiar} in connection with infinitesimal Hopf
algebras.

The AYBE is an analog in the world of associative algebras of the well-known classical Yang-Baxter equation (CYBE),
$$[r^{12}(v),r^{13}(v+v')]+[r^{12}(v),r^{23}(v')]+[r^{13}(v+v'),r^{23}(v')]=0,$$
where $r(v)$ is a meromorphic function in a neighborhood of $0$ taking values in a Lie algebra $\fg \otimes \fg$.
Solutions of the CYBE lead to Poisson-Lie groups and classical integrable systems (see for ex.
\cite{drinfeld}, \cite{charipressley}). There is a direct relation between the two equations in the
case $\fg=\ssl_n(\C)$: if $r(u,v)$ is a skew-symmetric solution of the AYBE such that the limit $\overline{r}(v) = (pr \otimes
pr)r(u,v)|_{u=0}$ exists (where $pr$ is the projection away from the identity to traceless
matrices), then $\overline{r}(v)$ is a solution of the CYBE for $\fg=\ssl_n(\C)$.

It was discovered in \cite{pol02} that solutions of the AYBE often arise from $1$-Calabi-Yau $A_\infty$-categories.
More precisely, assume we have such a minimal $A_\infty$-category $\CC$ and two sets of isomorphism classes
of objects 
in $\CC$, $\mathcal{X}$ and $\mathcal{Y}$,
such that for every pair of distinct objects $x_1,x_2 \in \mathcal{X}$ (resp. $y_1,y_2 \in \mathcal{Y}$), $\Hom^*(x_1,x_2)=0$ (resp. $\Hom^*(y_1,y_2)=0$). Furthermore, assume  
$\Hom^{\neq 0}(x,y)=0$ (and so $\Hom^{\neq 1}(y,x)=0$) for $x\in \mathcal{X}$ and $y\in \mathcal{Y}$. Then 
dualising the triple product \footnote{Our convention is that in an $A_\infty$-category, we read the
compositions from right-to-left (as in \cite{seidelbook}). This affects certain signs in
computations. In particular, the $A_\infty$-relations are given by:
\[ \sum_{m,n} (-1)^{|a_1|+\ldots + |a_n|-n} \m_{d-m+1}(a_d,\ldots, a_{n+m+1},      \m_m (a_{n+m},\ldots, a_{n+1}), a_n,\ldots a_1) = 0 \]
}
$$\m_3:\Hom^0(x_2,y_2) \ot\Hom^1(y_1,x_2)\ot\Hom^0(x_1,y_1)\to\Hom^0(x_1,y_2),$$
where $x_1,x_2\in \mathcal{X}$, $y_1,y_2\in \mathcal{Y}$, using the Calabi-Yau pairing, we get a tensor 
$$r^{x_1,x_2}_{y_1,y_2}:\Hom^0(x_2,y_2)\ot \Hom^0(x_1,y_1)\to\Hom^0(x_1,y_2)\ot \Hom^0(x_2,y_1).$$
defined by
\begin{equation} \label{dualiz} \lan r^{x_1,x_2}_{y_1,y_2}(f_{22} \otimes f_{11}), g_{21} \otimes
    g_{12} \ran 
= \lan \m_3(f_{22},g_{12}, f_{11}) , g_{21} \ran . \end{equation}
where $f_{ii} \in \Hom^0 (x_i,y_i), g_{ij} \in \Hom^1 (y_i,x_j)$.
Note that by the cyclicity of the $A_\infty$-structure, this tensor satisfies the following {\it skew-symmetry condition}:
\begin{equation}\label{skew-sym-eq}
(r^{x_1x_2}_{y_1y_2})^{21}=-r^{x_2x_1}_{y_2y_1}.
\end{equation}

Now, let $x_1,x_2,x_3$ (resp., $y_1,y_2,y_3$) be distinct elements of $\mathcal{X}$ (resp., $\mathcal{Y}$).
Then taking into account the assumptions on $\mathcal{X}, \mathcal{Y}$ the relevant
$A_\infty$-relation involving $\m_3$ takes the form 
\begin{equation}\label{m3-m3-id}
\m_3(\m_3(f_{33},g_{23},f_{22}),g_{12},f_{11})+\m_3(f_{33},\m_3(g_{23},f_{22},g_{12}),f_{11})-\m_3(f_{33},g_{23},\m_3(f_{22},
g_{12},f_{11}))=0,
\end{equation}
where $f_{ii}\in\Hom^0(x_i,y_i)$, $g_{ij}\in\Hom^1(y_i,x_j)$. 
Note that here the first and the the third terms can be immediately expressed in terms of the tensor $r^{x_1x_2}_{y_1y_2}$. 
To do this for the
middle term, one has to use the cyclic symmetry of the $A_\infty$-structure, which gives
$$\lan f_{31},\m_3(g_{23},f_{22},g_{12})\ran=\lan g_{12},\m_3(f_{31},g_{23},f_{22})\ran.$$

Taking into account the cyclic symmetry, we can rewrite the above $A_\infty$-relation as follows 
(see \cite[Thm. 1]{pol02}):
\begin{equation}\label{Rid0}
(r^{x_1x_2}_{y_1y_3})^{13}(r^{x_2x_3}_{y_2y_3})^{12}+
(r^{x_3x_2}_{y_1y_2})^{23}(r^{x_1x_3}_{y_1y_3})^{13}-
(r^{x_1x_3}_{y_2y_3})^{12}(r^{x_1x_2}_{y_1y_2})^{23}=0.
\end{equation}
This is viewed as an equation on 
\[ \Hom^0(x_3,y_3) \ot \Hom^0(x_2,y_2) \ot \Hom^0(x_1,y_1) \to  \Hom^0(x_2,y_3) \ot \Hom^0(x_1,y_2) \ot \Hom^0(x_3,y_1) \]
Permuting the second and third factors in the tensor product, and swapping
$x_1$ with $x_2$ and $y_1$ and $y_2$ (and also taking into account the skew-symmetry \eqref{skew-sym-eq}), 
the equation \eqref{Rid0} is equivalent to the following equation
\begin{equation}\label{Rid}
(r^{x_2x_1}_{y_2y_3})^{12}(r^{x_1x_3}_{y_1y_3})^{13}-
(r^{x_1x_3}_{y_1y_2})^{23}(r^{x_2x_3}_{y_2y_3})^{12}+
(r^{x_2x_3}_{y_1y_3})^{13}(r^{x_1x_2}_{y_1y_2})^{23}=0.
\end{equation}
We will call the equation \eqref{Rid} the {\it general AYBE} (or simply {\it AYBE} when no confusion can arise).\footnote{Our equation differs from 
\cite[Eq.\ (1.2)]{pol02} due to different conventions. 
The two equations become equivalent if we replace $r^{xx'}_{yy'}$ by $r^{x'x}_{y'y}$.}

It was further shown in \cite{pol02} that in the case when $\CC$ is the derived category of coherent sheaves on an elliptic
curve (or some of its degenerations) then there exist natural choices of $\mathcal{X}$ and $\mathcal{Y}$ as above, so that all the spaces
$\Hom^0(x,y)$, $x\in \mathcal{X}$, $y\in \mathcal{Y}$, can be identified with the fixed finite-dimensional vector space $V$.
Furthermore, in this case $\mathcal{X}$ and $\mathcal{Y}$ have abelian group structures, and
the obtained tensors $r^{x_1,x_2}_{y_1,y_2}:V^{\ot 2}\to V^{\ot 2}$ depend only on the differences
$u=x_2-x_1$, $v=y_2-y_1$,
which leads to the equation \eqref{AYBE}. 

Note that different choices of trivialization of the $\Hom$-spaces in the above 
construction correspond to the natural equivalence relation on solutions of the AYBE introduced in \cite{pol02}.
Namely, given a function $\varphi^x_y$ with values in $\GL_n(\C)$, we can transform a solution $r^{x_1x_2}_{y_1y_2}$
of \eqref{Rid} to the new solution
\begin{equation}\label{equivalence-eq}
\wt{r}^{x_1x_2}_{y_1y_2}=(\varphi^{x_2}_{y_1}\ot\varphi^{x_1}_{y_2})r^{x_1x_2}_{y_1y_2}
(\varphi^{x_1}_{y_1}\ot\varphi^{x_2}_{y_2})^{-1}.
\end{equation}

Our first result is that an analog of the above construction gives a bijection between formal
solutions of the general AYBE and a class of $A_\infty$-structures. Namely, we will consider
deformations of the formal $A_\infty$-category $\AA=\AA_n$ defined below. Note that we use the sign conventions
of \cite{seidelbook}, so that the double compositions in the associated cohomology category differ from 
those induced by $\m_2$ by a sign.

\begin{defi}\label{A-cat-defi}
The $A_\infty$-category $\AA=\AA_n$ has two objects $X$ and $Y$, and the $\Hom$-spaces
$$\Hom(X,Y)=\Hom^0(X,Y)=\Z\th_1\oplus\ldots\oplus\Z\th_n, \ \
\Hom(Y,X)=\Hom^1(Y,X)=\Z\eta_1\oplus\ldots\oplus\Z\eta_n,$$ $$\Hom^0(X,X)=\Z\id_X, \ \
\Hom^0(Y,Y)=\Z\id_Y, \ \ \Hom^1(X,X)=\Z\xi_X, \ \ \Hom^1(Y,Y)=\Z\xi_Y.$$ The elements $\id_X$ and
$\id_Y$ act as strict units in the sense that \[ \m_2(a , \id_X) = a \ \ , \ \ \m_2( \id_X, a) =
(-1)^{|a|} a \ \ , \ \ \m_2(a , \id_Y) = a \ \ , \ \  \m_2( \id_Y, a) = (-1)^{|a|} a, \] whenever
composition with $a \in \mathcal{A}$ is non-zero, where $|a|$ is the degree of $a$, and the other compositions are given by 
$$\m_2(\eta_\a, \th_\b)=\de_{\a\b}\xi_X, \ \ \m_2( \th_\a,
\eta_\b)=-\de_{\a\b}\xi_Y.$$
\end{defi}

Note that one can view $\mathcal{A}$ as a graded category by defining the composition as:
\[ a_2 \cdot a_1 = (-1)^{|a_1|} \m_2(a_2,a_1) \]
We also define the symmetric perfect pairing on the $\Hom$-spaces of $\AA$ by
$$\lan \eta_\a,\th_\b\ran= - \lan \th_b , \eta_a \ran = \de_{\a\b}, \ \ \lan \xi_X,\id_X\ran=- \lan
\id_X, \xi_X \ran = \lan \xi_Y,\id_Y\ran =-\lan \id_Y, \xi_Y \ran = 1.$$
Let $\k$ be a field. We are going to consider $A_\infty$-structures on
$\AA \otimes \k$, with given $\m_2$, which are cyclic with respect to this pairing. 
Recall that a {\it strictly cyclic $A_\infty$-category of dimension 1} is a strictly unital, proper
$A_\infty$-category together with nondegenerate pairings 
\[ \langle\ ,\ \rangle : hom^*(X,Y) \otimes hom^{1-*} (Y,X) \to \k \]
satisfying
\[ \langle a_1 , a_2 \rangle = (-1)^{(|a_1|-1)(|a_2|-1)+1} \langle a_2, a_1 \rangle \] 
and the cyclic symmetry condition:
\[ \langle a_{k+1}, \m_k (a_k, a_{k-1}, \ldots, a_1) \rangle =
(-1)^{(|a_{k+1}|-1|)(|a_1|+|a_2|+\ldots + |a_k|-k)} \langle a_k, \m_k (a_{k-1},a_{k-2}\ldots, a_1,
a_{k+1}) \rangle \]

An $A_\infty$-functor $\mathfrak{f}= (\mathfrak{f}^n)_{n\geq 1}: \mathscr{A} \to \mathscr{B}$
between cyclic $A_\infty$ categories is said to be \emph{cyclic} if the following hold:
\[ \langle \mathfrak{f}^1(a_2), \mathfrak{f}^1(a_1) \rangle = \langle a_2, a_1 \rangle \]
for any $a_2,a_1$ and 
\[ \sum_{k+l=n} \langle \mathfrak{f}^l (a_n,\ldots, a_{k+1}) , \mathfrak{f}^k (a_k,\ldots a_1)
\rangle =0 \] 
for any sequence of composable morphisms $a_n,\ldots, a_1$ with $n \geq 3$.

For a commutative $\k$-algebra $R$ we denote by $\MM_\infty(\AA\ot R)$ the set of
cyclic, strictly unital, minimal $A_\infty$-structures on $\AA\ot R$, up to a strict cyclic $A_\infty$-equivalence (i.e., the one
with $\mathfrak{f}^1=\id$).
Let us set 
\begin{equation}\label{P-def}
    \mathrm{P}:=\sum_{i,j} e_{ij}\ot e_{ji}\in \Mat_n(\k)\ot\Mat_n(\k),
\end{equation}
where $(e_{ij})$ is the standard basis of $\Mat_n(\k)$ defined by
$e_{ij}(\mathbf{e}_k) = \delta_{jk} \mathbf{e}_i$ if $(\mathbf{e}_i)_{i=1}^n$ is a basis of $\k^n$. In other words, $\mathrm{P}$ is the transposition
operator given by: \[ \mathrm{P} (x\otimes y) = y \otimes
x. \]

\noindent
{\bf Theorem A}. {\it There is a natural explicit bijection between $\MM_\infty(\AA\ot R)$ and the equivalence classes
of formal skew-symmetric
solutions $r^{x_1x_2}_{y_1y_2}$ of the general AYBE of the following type. We let $x_1,x_2,y_1,y_2$ to be formal
variables and consider 
$$r^{x_1x_2}_{y_1y_2}\in \Mat_n(\k)\ot\Mat_n(\k)\ot
R[[x_1,x_2,y_1,y_2]][(x_2-x_1)^{-1}(y_2-y_1)^{-1}]$$
of the form
\begin{equation}\label{r-expansion-eq}
    r^{x_1x_2}_{y_1y_2}\equiv \frac{\id\ot\id}{x_2-x_1}+\frac{\mathrm{P}}{y_1-y_2} \ \ \  \{ \mod \ \
    \Mat_n(\k)\ot\Mat_n(\k)\ot R[[x_1,x_2,y_1,y_2]] \} ,
\end{equation}
such that \eqref{Rid} is satisfied in 
$\Mat_n(\k)\ot\Mat_n(\k)\ot R[[x_1,x_2,x_3,y_1,y_2,y_3]][\De^{-1}]$, where
$\De=\prod_{i<j}(x_j-x_i)(y_j-y_i)$.
The skew-symmetry is the equation \eqref{skew-sym-eq}.
The equivalence between such solutions is given by \eqref{equivalence-eq}, with
$$\varphi^x_y\in \mathrm{Id} + (x,y)\sub \Mat_n(R)[[x,y]].$$
Considering the more general equivalences, where the constant term of $\varphi^x_y$ is only required to be an invertible matrix,
corresponds to general cyclic $A_\infty$-equivalences of the cyclic $A_\infty$-structures.
}

The key idea in this theorem is to apply a version of the above construction of solutions of the AYBE to a pair
$\XX$, $\YY$ of formal deformations of objects $X$ and $Y$ in $\AA$. A technical point is that these formal
deformations are defined in the category of twisted objects over $\AA$, which is non-minimal. Because of this
one has to use certain triple Massey products instead of just $\m_3$ (see Sec.\ \ref{Massey-sec}).
In particular, the singular terms in the expansion of $r^{x_1x_2}_{y_1y_2}$ are obtained naturally in this approach
due to the definition of the Massey products.

Recall that Belavin and Drinfeld in the seminal paper \cite{BD} classified 
nondegenerate\footnote{this means that the tensor $r(v)$ is nondegenerate for generic $v$}  solutions of the classical
Yang-Baxter equation for simple complex Lie algebras, up to some natural equivalence. They showed that they can be
 either elliptic or trigonometric or rational, and further classified trigonometric solutions in terms of some combinatorial data,
involving so called Belavin-Drinfeld triples. 

Similarly, one can pose the problem of classifying nondegenerate solutions $r(u,v)$ of the AYBE 
                                     (and of its formal general
version). Partial results in this direction we obtained in \cite{pol02} and \cite{pol09}. If
$r$ is strongly nondegenerate (see Def. \ref{strdegdef} and Prop. \ref{strdeg}), the Laurent expansion of the solution at $u=0$ has the form 
\begin{equation}\label{Laurent-exp-ass}
r(u,v)=\frac{1\ot 1}{u}+r_0(v)+\ldots
\end{equation}
Under this assumption, it was shown that the projection $\ov{r}_0(v)$
of $r_0(v)$ to $\ssl_n(\C)\ot\ssl_n(\C)$ is a nondegenerate solution of the CYBE, and that if $r_0(v)$ is either elliptic or trigonometric then $r(u,v)$ is determined by $\ov{r}_0(v)$, up to some natural transformations.
Note that the Laurent expansion \eqref{Laurent-exp-ass} appears naturally in the construction of Theorem A.
It was shown in \cite{pol02} that all elliptic solutions of the CYBE extend to those of the AYBE.
In \cite{sch} Schedler observed that this is not the case for all the trigonometric solutions. Extending this work, it was proved in \cite{pol09} that nondegenerate
solutions of the AYBE, with the Laurent expansion at $u=0$ of the form \eqref{Laurent-exp-ass} and such that
$\ov{r}_0(v)$ is a trigonometric solution of the CYBE, admit a classification in terms of the following combinatorial data
(see also Sec.\ \ref{Fukaya} below).

\begin{defi} \label{BD} An \emph{associative Belavin-Drinfeld structure} $(S,C_1,C_2,A)$ consists of a finite set $S$, a pair of transitive permutations $C_1, C_2 : S \to S$ and a proper subset $A \subset S$ such that for all $a \in A$, one has :
	\[ C_1(C_2(a)) = C_2(C_1(a)). \]
\end{defi} 

The reader familiar with the original Belavin-Drinfeld triples (defined in terms of Dynkin diagrams)
may notice that the above associative analog of this notion is more elementary (the original definition in \cite{pol09} 
is slightly different but
is equivalent to the one above, see Sec.\ \ref{Fukaya}). 

One can also ask which solutions of the AYBE can be realized by families $(\XX,\YY)$ of objects
in some geometric $1$-Calabi-Yau-categories. A natural source is provided by the derived categories of coherent sheaves
on elliptic curves and their degenerations. Then we can take as $\XX$ a universal deformation of a simple vector bundle,
and as $\YY$ the family of the structure sheaves of points.

It turns out that all the solutions of the AYBE for which $\ov{r}_0(v)$ is elliptic
arise in this way from simple vector bundles on elliptic curve, and can be explicitly computed in terms
of elliptic functions (see \cite{pol02}). In \cite{pol09}, all the solutions coming from
the nodal degenerations of elliptic curves, i.e., cycles of projective lines (aka standard $m$-gons), were computed
and were shown to be trigonometric. However, it turned out, somewhat unexpectedly, 
that not all trigonometric solutions of the AYBE appear in this way.
Namely, it was also shown in \cite{pol09} that the trigonometric
solution of the AYBE, corresponding to the data $(S,C_1,C_2,A)$, 
arises from a simple vector bundle on a
cycle of projective lines if and only if the corresponding transitive permutations $C_1$ and $C_2$ commute (equivalently,
$C_2=C_1^k$ for some $k$).

This raised a natural problem of finding other $1$-Calabi-Yau $A_\infty$-categories and objects in them,
which would account for missing solutions. This problem is solved in the present paper by looking at appropriate 
Fukaya categories. Namely, starting from the data of an associative Belavin-Drinfeld structure
$(S,C_1,C_2,A)$, we construct a square-tiled surface $\Sigma$ with a local
symplectomorphism
\[ \pi : \Sigma \to \mathbb{T} \]
to the square torus $\mathbb{T}$. In the case $A = \emptyset$, $\Sigma$ is just the $n-$fold
covering space of the punctured torus $\mathbb{T}_0$ associated to the permutations $C_1, C_2$ (see Section \ref{symplectic} for the general case). Lifts of standard Lagrangian curves in $\mathbb{T}$ to
$\Sigma$ give a pair of exact Lagrangians $L_1$ and $L_2$ in $\Sigma$ such that 
\[ \bigoplus_{i,j=1}^2 \mathrm{HF^*} (L_i,L_j) \simeq \mathcal{A} \otimes \C. \] 

Now, we have complex push-offs of the Lagrangians $L_1$ and $L_2$ forming
1-parameter families $L_1^{x}$ and $L_2^{y}$ (see Definition \ref{pushoff}). Taking these two families as families
                                     $\XX$ and $\YY$ in our general construction of solutions of the AYBE,
                                     we get such a solution that records
                                     triple products between
$(L_1^{x_1}, L_2^{y_1}, L_1^{x_2}, L_2^{y_2})$. We show that this gives
exactly the trigonometric solution of the AYBE associated with $(S,C_1,C_2,A)$. More precisely, we have:

{\bf Theorem B}. {\it Let $(\Sigma, L_1, L_2)$ be the square-tiled surface $\Sigma$ and Lagrangians $L_1$ and $L_2$ associated with an associative
Belavin-Drinfeld structure $(S,C_1,C_2,A)$. Then the tensor $r^{x_1,x_2}_{y_1,y_2}$ obtained from the triple products in the Fukaya category
$\mathcal{F}(\Sigma)$ only depends on $u=x_2-x_1$, $v=y_2-y_1$ and is a solution of the
AYBE over $\C$ given explicitly by the following formula (for an appropriate choice of basis):
\begin{equation}
\begin{aligned}
 r(u,v)&= \frac{1}{\exp(u)-1}\sum_{i} e_{ii}\ot e_{ii} + \frac{1}{1-\exp(-v)}\sum_i e_{ii}\ot e_{ii}  \\        
& +      \frac{1}{\exp(u)-1}\sum_{0< k<n, i}\exp(\frac{ku}{n})e_{C_1^k(i),C_1^k(i)}\ot e_{ii} +
    \frac{1}{\exp(v)-1}\sum_{0<m<n, i}\exp(\frac{mv}{n})e_{i,C_2^m(i)}\ot e_{C_2^m(i),i}  \\       &
    +\sum_{0<k,0<m; a\in A(k,m)}
\Big\{ \exp(-\frac{ku+mv}{n})e_{C_2^m(a),a}\ot
e_{C_1^k(a),C_1^k C_2^m(a)}-\exp(\frac{ku+mv}{n})e_{C_1^k(a),C_1^k C_2^m(a)}\ot e_{C_2^m(a),a}
\Big\}, \end{aligned}\end{equation}
    where 
    %$E(z) = e^z \in \C$ is the exponential function and 
    we denote by $A(k,m)\sub A$ the set of all $a\in A$ such that
$C_1^iC_2^j(a)\in A$ for all $0\le i< k, 0\le j< m$. }

We note that the surface $\Sigma$ has genus $1$ (i.e., is a punctured torus) if and only if $C_1$ and $C_2$
    commute. This explains why only these solutions appeared from simple vector bundles on nodal
    degenerations of elliptic curves, which are mirror dual to punctured tori (see e.g., \cite{LP}).

As an application of the viewpoint developed in this paper (combined with the results of \cite{pol09})
    we derive the following result about simple vector bundles on cycles of projective lines.

\medskip

\noindent
{\bf Theorem C}. {\it Let $C$ be a cycle of projective lines (the standard $n$-gon) over $\C$. For any simple vector bundle
$V$ on $C$ there exists a composition $\Phi$ of $1$-spherical twists and their inverses
                                                    such that $\Phi(\OO_C)\simeq V$.}

Here we use the notion of the twist autoequivalence associated with an $n$-spherical object 
                                introduced in \cite{Seidel-Thomas}. Recall that an $n$-spherical object $E$
                                                    should satisfy $\Hom^*(E,E)=k\oplus k[-n]$ (together with an
                                                    additional nondegeneracy condition). The corresponding
                                                    twist autoequivalence $T_E$ fits into an exact triangle
                                                    $$\Hom^*(E,X)\otimes E\to X\to T_E(X)\to\ldots$$
                                                    In this paper we consider only $1$-spherical objects and
                                                    the corresponding twists.
                                                    
Note that Theorem C is known in the case $n=1$ by the work \cite{BK}. In this case the situation is very similar to the case
of elliptic curves. The case $n>1$ is much more complicated: in this case one can still classify all simple vector
bundles on $C$ (see \cite{BDG}) but the relevant combinatorics is quite involved.

The idea of the proof of Theorem C is to consider the solution of the formal general AYBE associated with
the pair $(\OO,V)$ (where $V$ is sufficiently positive), and to use Theorem A which states that the subcategory
generated by $\OO$ and $V$ is determined by this solution. The point is that we know this solution of the AYBE to be the same
as for the Lagrangians $L_1,L_2$ in the symplectic surface $\Sigma$
of genus $1$ associated with some associative Belavin-Drinfeld structure.
We prove that in this situation the pair $(\OO,V)$ (resp. $(L_1,L_2)$) split generates the perfect derived category of $C$
(resp., the Fukaya category of $\Sigma$). Thus, we reduce the problem to a similar question about Lagrangians in
the Fukaya category, where we can use the action of the mapping class group.

The paper is organized as follows. In Section \ref{formal-AYBE-sec} we study the relation between
 formal solutions of the general AYBE and $A_\infty$-structures, in particular, proving Theorem A in \ref{proof-A-sec}.
 In addition, in Sec.\ \ref{involution-sec} we discuss the natural involution on solutions of the AYBE, which allows
 us to deduce the pole conditions imposed in \cite{pol09} for strongly nondegenerate solutions of the AYBE 
 (see Prop.\ \ref{strdeg}). Also, in Sec.\ \ref{alg-form-AYBE-sec} we explain, basing on ideas of \cite{Seidel-flux},
 the connection between solutions of the AYBE coming from algebraic (or analytic) families of
 objects (see \eqref{dualiz}) and the corresponding formal solutions from Theorem A.
Section \ref{symplectic} is devoted to the construction of trigonometric solutions of the AYBE from Fukaya products on
the square-tiled surfaces associated with Belavin-Drinfeld structures.
In Section \ref{cycles-sec} we consider two applications of Theorems A and B to vector bundles over a (nodal) cycle $C$
of projective lines. One is a criterion, in terms of some combinatorial data, for a pair of simple vector bundles on $C$
to be related by a Fourier-Mukai autoequivalence (see Theorem \ref{ABD-equiv-thm}). 
Another is Theorem C, proved in Sec.\ \ref{proof-C-sec}.

    {\it Acknowledgments.} Y.L. is supported in part by the Royal Society and the NSF grant
    DMS-1509141, and would like to thank Denis Auroux for a helpful correspondence.  A.P. is
    supported in part by the NSF grant DMS-1400390, and would like to thank SISSA, where part of this work was done,
    for hospitality and excellent working conditions. We are grateful to the referee for comments and suggestions.

\section{A class of cyclic $A_\infty$-structures and formal solutions of the general AYBE}\label{formal-AYBE-sec}

\subsection{Twisted objects over complete rings}
    \label{twistedobjects}

Let us quickly review the definition of the $A_\infty$-category of twisted objects with coefficients in a complete 
    ring (see \cite{FOOO-01}, \cite{Fuk-01}, \cite{Lefevre}), mostly following \cite[Ch. 1]{seidelbook} and \cite[Sec.\ 7.6]{Keller}.

Let $\CC$ be a topological 
$A_{\infty}$-category over a complete ring $R$. We assume that $R$ is topologized by a decreasing
    filtration
$(R_n)$, such that $R_mR_n\sub R_{m+n}$, and that the $\Hom$-spaces in $\CC$ are complete.
We will only consider twisted objects of the following kind: $(X,\de_X)$, where $X$ is an object of $\CC$ and 
$\de_X\in R_1\Hom^1(X,X)$ is an element satisfying the Maurer-Cartan equation
$$\sum_{n\ge 1} \m_n(\de_X^n)=0.$$
Note that here the left-hand side converges in $\Hom^2(X,X)$.
    The $\Hom$-space between two such objects $(X,\de_X)$ and $(Y,\de_Y)$ is simply $\Hom(X,Y)$.  
There are natural $A_\infty$-products $(\m_n^t)$ for the twisted objects, which are obtained by inserting the twisting elements
$\de$ any number of times wherever possible. More precisely, $\m_d^t$ is given by
$$\m_d^t(a_d,a_{d-1},\ldots, a_1)=\sum_{i_0,\ldots,i_d\ge 0}
    \m_{d+i_0+\ldots+i_d}(\de_{X_d}^{i_d},a_d,\de_{X_{d-1}}^{i_{d-1}},a_{d-1},\ldots,
    a_1,\de_{X_0}^{i_0})$$ (We follow the sign conventions of
    \cite[Ch.1]{seidelbook}).

Let us point out one additional feature of the $A_\infty$-category of twisted objects: it easy to check
  that if we start with a cyclic $A_\infty$-category then the corresponding $A_\infty$-category inherits the cyclic
structure.

An $A_\infty$-functor $\mathfrak{f} = (\mathfrak{f}^n)_{n\geq 1} :\CC\to \CC'$ between $A_\infty$-categories over $R$, as above,
leads to an $A_\infty$-functor $\mathcal{F}= (\mathcal{F}^n)_{n\geq 1}$ between their categories of twisted objects. Namely,
$(X,\de_X)$ maps to
$(\mathcal{F}(X),\mathcal{F}(\de_X))$, where
$$\mathcal{F}(X)= \mathfrak{f}(X) \ \ , \ \ \mathcal{F}(\de_X)=\sum \mathfrak{f}^n(\de_X^n).$$
The maps $\mathcal{F}^d$ on $\Hom$-spaces between twisted objects is given by
$$\mathcal{F}^d(a_d,a_{d-1},\ldots, a_1)=\sum_{i_0,\ldots,i_d\ge 0}
    \mathfrak{f}^{d+i_0+\ldots+i_d}(\de_{X_d}^{i_d},a_d,\de_{X_{d-1}}^{i_{d-1}},a_{d-1},\ldots,
    a_1,\de_{X_0}^{i_0})$$

\subsection{Triple Massey products and a construction of formal solutions of the general AYBE}
\label{Massey-sec}

First, let us recall the general definition of the triple Massey product in an $A_\infty$-category
$\CC$ following 
\cite{pol02}, but with different sign conventions.
For a triple of $\m_1$-closed composable morphisms 
$$X_0\rTo{a_1}X_1\rTo{a_2}X_2\rTo{a_3} X_3$$ 
such that $\m_2(a_2,a_1)=\m_1(h_1)$, $\m_2(a_3,a_2)=\m_1(h_2)$, one sets
\begin{equation}\label{MP-eq}
    \mathrm{MP}(a_3,a_2,a_1)=\m_3(a_3,a_2,a_1) - \m_2(h_2,a_1)- \m_2(a_3,h_1) \ \ \ \  \{ \mod \mathrm{Im}(\m_1) \},
\end{equation}
which is well-defined as a coset for the image of 
$$H^{|a_1|+|a_2|-1}\Hom(X_0,X_2)\oplus
H^{|a_2|+|a_3|-1}\Hom(X_1,X_3)\rTo{(\m_2(a_3,?),\m_2(?,a_1))}
H^{|a_1|+|a_2|+|a_3|-1}\Hom(X_0,X_3).$$

  The main result about such triple Massey products is that they are preserved under $A_\infty$-functors:
if $\mathfrak{f}:\CC\to \CC'$ is an $A_\infty$-functor then 
$$\mathfrak{f}^1(\mathrm{MP}(a_3,a_2,a_1))=\mathrm{MP}(\mathfrak{f}^1(a_3),\mathfrak{f}^1(a_2),\mathfrak{f}^1(a_1))$$
in the appropriate quotient-group (see \cite[Prop.\ 1.1]{pol02}).
In particular, if $\CC'$ is a minimal model of $\CC$ obtained as a result of the homological perturbation procedure,
  then the Massey product $MP(a_3,a_2,a_1)$, computed in $\CC$, agrees with $\m_3^{\CC'}(a_3,a_2,a_1)$
  (since for a minimal $A_\infty$-category our Massey product reduces to $\m_3$).

Thus, the construction of solutions of the AYBE from two families of objects, presented in
  the introduction, has a version for non-minimal cyclic $A_\infty$-categories, linear over some commutative ring $R$.
  More precisely, we have to replace in this construction $\m_3$ by the triple Massey product $\mathrm{MP}$ and
  assume that cohomology of all the morphism spaces are projective $R$-modules (so that the homological perturbation 
  can be applied).
 One technical problem is that the cyclic property of the $A_\infty$-structure is not necessarily inherited by the minimal model.
However, we have the following compatibility of the Massey products with the cyclic structure.

\begin{lem}\label{MP-cyclic-lem} 
Suppose we are given a cycle of $\m_1$-closed morphisms in a cyclic $A_\infty$-category,
$$X_0\rTo{a_1}X_1\rTo{a_2}X_2\rTo{a_3} X_3\rTo{a_4} X_0$$ 
such that 
$$\m_2(a_2,a_1)=\m_1(h_1),  \ \ \m_2(a_3,a_2)=\m_1(h_2), \ \ \m_2(a_4,a_3)=\m_1(h_3).$$
Assume also that the corresponding Massey products $\mathrm{MP}(a_3,a_2,a_1)$ and $\mathrm{MP}(a_4,a_3,a_2)$ are
univalued. Then
$$\lan a_1,\mathrm{MP}(a_4,a_3,a_2)\ran=(-1)^{(|a_1|-1)(|a_2|+|a_3|+|a_4|-1)}\lan a_4, \mathrm{MP}(a_3,a_2,a_1)\ran.$$
\end{lem}

\Pf . Using \eqref{MP-eq}, we see that it is enough to establish the following identities
$$\lan a_1,\m_3(a_4,a_3,a_2)\ran=(-1)^{(|a_1|-1)(|a_2|+|a_3|+|a_4|-1)}\lan a_4,\m_3(a_3,a_2,a_1)\ran,$$
$$\lan a_1,\m_2(a_4,h_2)\ran=(-1)^{(|a_1|-1)(|a_2|+|a_3|+|a_4|-1)}\lan a_4,\m_2(h_2,a_1)\ran,$$
$$\lan a_1,\m_2(h_3,a_2)\ran=(-1)^{(|a_1|-1)(|a_2|+|a_3|+|a_4|-1)}\lan a_4,\m_2(a_3,h_1)\ran.$$
The first two follow directly from the cyclic property of $\m_3$ and $\m_2$ (noting that $|h_2|=|a_2|+|a_3|-1$).
For the last one, we first rewrite the left-hand side using the cyclic property of $\m_2$:
$$\lan a_1,\m_2(h_3,a_2)\ran=(-1)^{(|a_1|-1)(|a_2|+|a_3|+|a_4|-1)}\lan h_3,\m_2(a_2,a_1)\ran=
(-1)^{(|a_1|-1)(|a_2|+|a_3|+|a_4|-1)}\lan h_3,\m_1(h_1)\ran.$$
Next, we use the cyclic property of $\m_1$:
$$\lan h_3,\m_1(h_1)\ran=(-1)^{(|h_1|-1)(|h_3|-1)}\lan h_1,\m_1(h_3)\ran=
(-1)^{(|h_1|-1)(|h_3|-1)}\lan h_1,\m_2(a_4,a_3)\ran.$$
It remains to apply the cyclic of $\m_2$ to the last expression to get the required identity.
\ed
  
Using Lemma \ref{MP-cyclic-lem}, one can generalize the construction of the solutions of the AYBE to non-minimal
  cyclic $A_\infty$-categories provided the appropriate vanishing assumptions hold on the cohomology level.
  
We want to apply this construction to certain Massey products involving twisted objects over an
  $A_\infty$-structure on the category $\AA$ (see Def.\ \ref{A-cat-defi})
 More precisely, given a minimal cyclic $A_\infty$-structure with the given $m_2$ 
on the category $\AA\ot R$, where $R$ is a commutative ring,
we extend the coefficients to 
$$\wt{R}:=R[[x_1,x_2,y_1,y_2]][(x_2-x_1)^{-1}(y_2-y_1)^{-1}]$$
and consider the following twisted objects of $\AA\ot\wt{R}$:
\begin{equation}\label{Xi-Yi-def}
X_i=(X,x_i\xi_X), \ \ Y_i=(Y,y_i\xi_Y), \ \text{ for } i=1,2.
\end{equation}

Note that the $\Hom$-spaces $\Hom(X_i,Y_j)$ and $\Hom(Y_j,X_i)$ are still concentrated in one degree
and so have trivial $\m_1^t$. We denote by $\th_\a^{ij}$ (resp., $\eta_\a^{ji}$) the basis elements in $\Hom(X,Y)$ (resp,
$\Hom(Y,X)$) viewed as elements of $\Hom(X_i,Y_j)$ (resp., $\Hom(Y_j,X_i)$).
On the other hand, $\Hom(X_1,X_2)$ and $\Hom(Y_1,Y_2)$ now have a nontrivial differential:
$$\m_1^t(\id_X)=(x_2-x_1)\xi_X, \ \ \m_1^t(\id_Y)=(y_2-y_1)\xi_Y,$$
so the corresponding cohomology vanishes (due to the localization in the definition of $\wt{R}$).

We consider the triple Massey product corresponding to
the composable morphisms
$$X_1\rTo{\th_\a^{11}}Y_1\rTo{\eta_\b^{12}}X_2\rTo{\th_{\a'}^{22}}Y_2.$$
We claim that this Massey product is well-defined and univalued. Indeed, we have
$$\m_2^t(\eta_\b^{12},\th_\a^{11})=\de_{\a\b}\xi_X=\frac{\de_{\a\b}}{x_2-x_1}\cdot \m_1^t(\id_X),$$
$$\m_2^t(\th_{\a'}^{22},\eta_\b^{12})=- \de_{\a'\b}\xi_Y= \frac{\de_{\a'\b}}{y_1-y_2}\cdot \m_1^t(\id_Y),$$
hence it is well-defined. The fact that it is univalued follows immediately from the vanishing
of $H^0\Hom(X_1,X_2)$ and $H^0\Hom(Y_1,Y_2)$.
According to the formula \eqref{MP-eq} we have
\begin{align*}
    \mathrm{MP}(\th_{\a'}^{22},\eta_\b^{12},\th_{\a}^{11}) &=
    \m_3^t(\th_{\a'}^{22},\eta_\b^{12},\th_\a^{11}) -\frac{\de_{\a\b}}{x_2-x_1}\cdot
    \m_2^t(\th_{\a'}^{22}, \id_X) - \frac{\de_{\a'\b}}{y_1-y_2}\cdot \m_2^t(\id_Y, \th_\a^{11}) \\
    &= \m_3^t(\th_{\a'}^{22},\eta_\b^{12},\th_\a^{11}) -\frac{\de_{\a\b}}{x_2-x_1}\cdot
    \th_{\a'}^{22} -\frac{\de_{\a'\b}}{y_1-y_2}\cdot \th_\a^{11},  
\end{align*}
where the last equality follows since our $A_\infty$-structure on $\AA\ot R$ is strictly unital,
so the products $\m_2^t$ involving the identity remain the same as  $\m_2$. Therefore, we
have:
\begin{equation}\label{twisted-MP-eq}
    \lan \mathrm{MP}(\th_{\a'}^{22},\eta_\b^{12},\th_{\a}^{11}),\eta_{\b'}^{21}\ran=
    \lan \m_3^t(\th_{\a'}^{22},\eta_\b^{12},\th_{\a}^{11}),\eta_{\b'}^{21}\ran
+\frac{\de_{\a\b}\de_{\a'\b'}}{x_2-x_1}+
\frac{\de_{\a'\b}\de_{\a\b'}}{y_1-y_2}.
\end{equation}

Since the $A_\infty$-category of twisted objects is still cyclic, 
one can show that the above triple Massey product gives rise to a solution of the general AYBE  
over $R[[x_1,x_2,x_3,y_1,y_2,y_3]][\De^{-1}]$, which would prove one part of Theorem A.
Namely, one first shows that an analog of the $A_\infty$-identity \eqref{m3-m3-id}
  holds for the Massey products by passing to an equivalent minimal $A_\infty$-structure, and then
  uses Lemma \ref{MP-cyclic-lem} to rewrite the middle term. Even though this construction was what led us
  to Theorem A, we will use a different argument in its proof (since we need to show both directions).

\subsection{Proof of Theorem A}\label{proof-A-sec}
  To a minimal cyclic $A_\infty$-structure on $\AA\ot R$ we associate the element
$r^{x_1x_2}_{y_1y_2}\in \Mat_n(\k)\ot\Mat_n(\k)\ot \wt{R}$ obtained from 
$\mathrm{MP}(\th_{\a'}^{22},\eta_\beta^{12},\th_{\a}^{11})$ by dualization. In other words,
\begin{equation}\label{ainf-r-eq}
    r^{x_1x_2}_{y_1y_2}=\sum_{\a,\a',\b,\b'} \lan
    \mathrm{MP}(\th_{\a'}^{22},\eta_\b^{12},\th_{\a}^{11}),\eta_{\b'}^{21}\ran\cdot
e_{\b'\a'}\ot e_{\b\a}.
\end{equation}

Let us set
$$f^{x_1x_2}_{y_1y_2}:=
\sum_{\a,\a',\b,\b'} \lan \m_3^t(\th_{\a'}^{22},\eta_\b^{12},\th_{\a}^{11}),\eta_{\b'}^{21}\ran\cdot
e_{\b'\a'}\ot e_{\b\a},$$
which is an element of $ \Mat_n(\k)\ot\Mat_n(\k)\ot R[[x_1,x_2,y_1,y_2]]$.
We can rewrite \eqref{twisted-MP-eq} as
\begin{align*}
&r^{x_1x_2}_{y_1y_2}=f^{x_1x_2}_{y_1y_2}+
\frac{1}{x_2-x_1}\cdot\sum_{\a,\a'}e_{\a\a}\ot e_{\a'\a'}+\frac{1}{y_1-y_2}\cdot\sum_{\a,\a'}e_{\a'\a}\ot e_{\a\a'}=\\
&f^{x_1x_2}_{y_1y_2}+\frac{\id\ot\id}{x_2-x_1}+\frac{P}{y_1-y_2},
\end{align*}
In particular, the singular part of $r^{x_1x_2}_{y_1y_2}$ has the required form. Also, the skew-symmetry equation
\eqref{skew-sym-eq} is equivalent
to 
$$(f^{x_1x_2}_{y_1y_2})^{21}=-f^{x_2x_1}_{y_2y_1},$$
which can be deduced from the cyclic symmetry equation as follows:
\begin{align*}
    (f^{x_1x_2}_{y_1y_2})^{21} &= \sum_{\alpha, \alpha', \beta,\beta'} \lan
    \m_3^t(\th_{\a'}^{22},\eta_\b^{12},\th_{\a}^{11}),\eta_{\b'}^{21}\ran\cdot
    e_{\b \a }\ot e_{\b'\a'} \\ &= -\sum_{\alpha, \alpha', \beta,\beta'} \lan
    \m_3^t(\th_{\a}^{22},\eta_{\b'}^{12},\th_{\a'}^{11}),\eta_{\b}^{21}\ran\cdot
e_{\b \a }\ot e_{\b' \a'} \\
   &= -\sum_{\alpha, \alpha', \beta,\beta'} \lan
    \m_3^t(\th_{\a'}^{22},\eta_{\b}^{12},\th_{\a}^{11}),\eta_{\b'}^{21}\ran\cdot
e_{\b' \a' }\ot e_{\b \a} = - f^{x_2 x_1}_{y_2 y_1}.
\end{align*}
We claim that the general AYBE for $r^{x_1x_2}_{y_1y_2}$ is equivalent to the $A_\infty$-constraints in
$\AA\ot R$ applied to 
all possible strings of composable elements 
\begin{equation}\label{string-eq}
\xi_{Y_1}^f,\eta^{31}_{\a_3},\xi_{X_3}^e,\th^{33}_{\b_2},\xi_{Y_3}^d,\eta^{23}_{\a_2},\xi_{X_2}^c,\th^{22}_{\b_1},\xi_{Y_2}^b,\eta^{12}_{\a_1},\xi_{X_1}^a. \end{equation}
Note that because of the cyclic symmetry these constraints are equivalent to the full $A_\infty$-constraints.

Recall that our general AYBE takes place over the ring $R[[x_1,x_2,x_3,y_1,y_2,y_3]][\De^{-1}]$.
Over this ring we can define twisted objects $(X_i,Y_i)$ for $i=1,2,3$ as in \eqref{Xi-Yi-def}.
We extend the notation $\th_\a^{ij}$ and $\eta_\b^{ji}$ for basis elements in the $\Hom$-spaces to this case.
Let us set for brevity
$$\la_{ij}=\frac{\id\ot\id}{x_j-x_i}, \ \ \mu_{ij}=\frac{P}{y_i-y_j}.$$
Thus, the AYBE takes the following form:
\begin{align*}
&(f^{x_2x_1}_{y_2y_3}+\la_{21}+\mu_{23})^{12}(f^{x_1 x_3}_{y_1y_3}+\la_{13} + \mu_{13})^{13} -
    (f^{x_1 x_3}_{y_1y_2} + \la_{13} +\mu_{12})^{23}(f^{x_2x_3}_{y_2y_3}+\la_{23} +\mu_{23})^{12}
\\
&+(f^{x_2 x_3}_{y_1y_3}+\la_{23}+\mu_{13})^{13}(f^{x_1x_2}_{y_1y_2}+\la_{12}+\mu_{12})^{23}=0.
\end{align*}
A straightforward calculation shows that the terms depending quadratically on $(\la_{ij},\mu_{ij})$ cancel out, so this is equivalent to
an equation of the form
\begin{equation}\label{AYBE-f-eq}
    \mathrm{AYBE}[f]+
    \left[(f^{x_2
    x_1}_{y_2y_3})^{12}(\la_{13})^{13}-(\la_{13})^{23}(f^{x_2x_3}_{y_2y_3})^{12}\right]+\ldots=0,
\end{equation}
where $\mathrm{AYBE}[f]$ is the left-hand side of \eqref{Rid} with $f$ instead of $r$, and the remaining terms similarly combine
the terms linear in $\la_{ij}$ and $\mu_{ij}$.

Now we claim that looking at the coefficients of the expansion of \eqref{AYBE-f-eq} in $x_1,x_2,x_3,y_1,y_2,y_3$ 
we get precisely the $A_\infty$-constraints for the strings \eqref{string-eq}. 
%The coefficient of $x_1^a y_2^b x_2^c y_3^d x_3^e y_1^f$ in 
These constraints have the form
\begin{align*}
-&\sum_{f=f_2+f_1;c=c_2+c_1}
    \m_*(\xi_{Y_1}^{f_2},\m_*(\xi_{Y_1}^{f_1},\eta_{\a_3},\xi_{X_3}^e,\th_{\b_2},\xi_{Y_3}^d,
    \eta_{\a_2},\xi_{X_2}^{c_2}),\xi_{X_2}^{c_1},\th_{\b_1},\xi_{Y_2}^b,\eta_{\a_1},\xi_{X_1}^a) \\
    - &\sum_{e=e_2+e_1;b=b_2+b_1}
    \m_*(\xi_{Y_1}^f,\eta_{\a_3},\xi_{X_3}^{e_2},\m_*(\xi_{X_3}^{e_1},\th_{\b_2},\xi_{Y_3}^d,\eta_{\a_2},\xi_{X_2}^c,\th_{\b_1},\xi_{Y_2}^{b_2}),\xi_{Y_2}^{b_1},\eta_{\a_1},\xi_{X_1}^a) \\ 
    + &\sum_{d=d_2+d_1;a=a_2+a_1}
    \m_*(\xi_{Y_1}^f,\eta_{\a_3},\xi_{X_3}^e,\th_{\b_2},\xi_{Y_3}^{d_2},\m_*(\xi_{Y_3}^{d_1},\eta_{\a_2},\xi_{X_2}^c,\th_{\b_1},\xi_{Y_2}^b,\eta_{\a_1},\xi_{X_1}^{a_2}),\xi_{X_1}^{a_1})
    \\ + 
    &\ldots =0  ,
\end{align*}
where the additional terms
appear when one of $a,b,c,d,e,f$ is zero, and have the form either $\m_*(\ldots,\m_2,\ldots)$ or $\m_2(\ldots,\m_*,\ldots)$.
Using cyclic symmetry, one can immediately check that the coefficients of $x_1^a
y_2^bx_2^cy_3^dx_3^ey_1^f$ in
the three terms in $-\mathrm{AYBE}[f]$ match the first three terms 
in the $A_\infty$-constraint above. 
%Since these first three terms can be written as
%\[ -\m_3^t( \m_3^t(\eta_{\a_3}, \th_{\b_2}, \eta_{\a_2}), \th_{\b_1},\eta_{\a_1}) - \m_3^t(
%\eta_{\a_3}, \m_3^t(\th_{\b_2}, \eta_{\a_2}, \th_{\b_1}),\eta_{\a_1})
%+\m_3^t(\eta_{\a_3}, \th_{\b_2}, \m_3^t(\eta_{\a_2}, \th_{\b_1},\eta_{\a_1})) \] 

Now let us show how the second term in \eqref{AYBE-f-eq} matches
the terms with $f=0$ in the $A_\infty$-constraint. The matching of the other terms is done similarly.
First, we observe that 
$$(f^{x_2 x_1}_{y_2y_3})^{12}(\la_{13})^{13}-(\la_{13})^{23}(f^{x_2x_3}_{y_2y_3})^{12}=
\frac{1}{x_3-x_1}(f^{x_2x_1}_{y_2y_3}-f^{x_2 x_3}_{y_2y_3})^{12}.$$
If we expand $f^{x_2x_1}_{y_2y_3}$ in powers of $x_1$,
$$f^{x_2x_1}_{y_2y_3}=\sum_{n\ge 0} f^{x_2,n}_{y_2y_3}\cdot x_1^n,$$
then the above expression becomes
$$\sum_{n\ge 0} \frac{x_1^n-x_3^n}{x_3-x_1}\cdot (f^{x_2,n}_{y_2y_3})^{12}=
-\sum_{a,e\ge 0}x_1^a x_3^e\cdot (f^{x_2, a+e+1}_{y_2y_3})^{12}.$$
Now one can easily check that this matches the contribution of the terms of the form
$$\m_2(\eta_{\a_3},\m_*(\xi_{X_3}^e,\th_{\b_2},\xi_{Y_3}^d,
\eta_{\a_2},\xi_{X_2}^{c},\th_{\b_1},\xi_{Y_2}^b,\eta_{\a_1},\xi_{X_1}^a)) $$
in the $A_\infty$-constraint.

Next, we claim that changing an $A_\infty$-structure by a cyclic homotopy transforms the corresponding
solution of general AYBE to an equivalent one, as in \eqref{equivalence-eq}, and that all equivalences appear in
this way. Namely, a cyclic homotopy $(\mathfrak{f}^n)$ gives rise to the formal series
$$\varphi^x_y:=\sum_{a,b\ge 0} \lan \mathfrak{f}^{a+b+1}(\xi_Y^b,\eta_\a,\xi_X^a),\th_\b \ran\cdot
x^a y^b\cdot e_{\b\a}$$
in $\Mat_n(R)[[x,y]]$. Note that since $\mathfrak{f}^1$ is the identity, the constant term of $\varphi^x_y$ is equal to $\id\in \Mat_n(R)$.
One can easily check that if $(\m'_n)$ is obtained from $(\m_n)$ by a cyclic homotopy
$(\mathfrak{f}_n)$ then the corresponding solutions of the general AYBE are related by \eqref{equivalence-eq}.

Finally, we claim that a cyclic $A_\infty$-equivalence $(\mathfrak{f}^n)$ is uniquely determined by the corresponding formal series $\varphi^x_y$,
which could be an arbitrary series with the constant term equal to the identity. Indeed, recall that the condition
for a strict $A_\infty$-functor to be cyclic
is that
\begin{equation}\label{cyclic-homotopy-eq}
    \sum_{k+l=n}  \lan \mathfrak{f}^l (a_n,\ldots,a_{k+1}),\mathfrak{f}^k (a_{k},\ldots,a_1)\ran=0
\end{equation}
for any sequence of composable morphisms $a_1,\ldots,a_n$, where $n\ge 3$.
For our $A_\infty$-category, the only potentially non-trivial values of $\mathfrak{f}^*$ are
$$\mathfrak{f}^* (\xi_Y^n,\eta_\a,\xi_X^m), \ \ \mathfrak{f}^* (\xi_X^n,\th_\b,\xi_Y^m), \ \
\mathfrak{f}^* (\xi_Y^p,\eta_\a,\xi_X^n,\th_\b,\xi_Y^m), \ \ 
\mathfrak{f}^* (\xi_X^p,\th_\b,\xi_Y^n,\eta_\a,\xi_X^m).$$
The constraints between them are given by \eqref{cyclic-homotopy-eq} applied to the following two kinds of composable strings
$$(\xi_X^c,\th_\b,\xi_Y^b,\eta_\a,\xi_X^a) \ \ \text{ and } \ \ 
(\xi_Y^c,\eta_\a,\xi_X^b,\th_\b,\xi_Y^a),$$
where $a+b+c>0$.
The first of these strings gives the identity
\begin{align*}
    &\lan \xi, \mathfrak{f}^* (\xi_X^{c-1},\th_\b,\xi_Y^b,\eta_\a,\xi_X^a)\ran +
    \lan \mathfrak{f}^* (\xi_X^c,\th_\b,\xi_Y^b,\eta_\a,\xi_X^{a-1}),\xi \ran \\ 
    &- \sum_{b_1+b_2=b,\ga} \lan \mathfrak{f}^*(\xi_X^{c},\th_\b,\xi_Y^{b_2}),\eta_\ga \ran \cdot
  \lan  f^*(\xi_Y^{b_1},\eta_\a,\xi_X^{a}),\th_\ga \ran\ = 0,
\end{align*}
where the first (resp., second) term appears only for $a\ge 1$ (resp., $c\ge 1$). We can rewrite these identities in
terms of the generating series
$$\wt{\varphi}^x_y:=\sum_{b,c\ge 0} \lan \mathfrak{f}^{b+c+1}(\xi_X^c,\th_\b,\xi_Y^b),\eta_\a
\ran\cdot x^c y^b\cdot e_{\b\a},$$
$$\psi^{x_1,x_2}_y:=\sum_{a,b,c \ge 0}\lan
\mathfrak{f}^{a+b+c+2}(\xi_X^c,\th_\b,\xi_Y^b,\eta_\a,\xi_X^a),\xi_X\ran\cdot
x_1^ay^bx_2^c\cdot e_{\b\a},$$
as follows:
$$(x_1-x_2)\psi^{x_1,x_2}_y-\wt{\varphi}^{x_2}_y\cdot \varphi^{x_1}_y -\mathrm{Id}=0.$$
Note that subtracting $\mathrm{Id}$ here corresponds to avoiding the case $a=b=c=0$ in cyclic homotopy equation.
Setting $x_2=x_1$ we deduce from this that
\begin{equation}\label{cyclic-homotopy-sol-1}
\wt{\varphi}^x_y=-(\varphi^x_y)^{-1},
\end{equation}
and so the above identity can be solved for $\psi^{x_1,x_2}_y$:
$$\psi^{x_1,x_2}_y=\frac{\mathrm{Id}-\varphi^{x_1}_y(\varphi^{x_2}_y)^{-1}}{x_1-x_2}.$$
Similarly, the constraints associated with the strings $(\xi_Y^c,\eta_\a,\xi_X^b,\th_\b,\xi_Y^a)$ boil down to
\eqref{cyclic-homotopy-sol-1} and to an equation expressing all
$\mathfrak{f}^*(\xi_Y^c,\eta_\a,\xi_X^b,\th_\b,\xi_Y^a)$ in terms
of the coefficients of $\varphi^x_y$.
\ed

\subsection{Involution on skew-symmetric solutions of the AYBE}\label{involution-sec}

In the following Proposition we define a natural involution on the skew-symmetric solutions of the general AYBE
\eqref{Rid}, where the variables $x_i$ and $y_i$ could be either distinct elements of some sets $\XX$ and $\YY$, or
formal variables, as in Theorem A.

\begin{prop}\label{propsymm} 
    Let $r^{x_1x_2}_{y_1y_2}$ be a solution of the general AYBE with values in $\mathrm{Mat}_n(\k) \otimes \mathrm{Mat}_n (\k) \otimes R$, where $R$ is a commutative $\k$-algebra,
satisfying the skew-symmetry condition \eqref{skew-sym-eq}. Set
$$\hat{r}^{\hat{x}_1\hat{x}_2}_{\hat{y}_1\hat{y}_2}:=(r^{\hat{y}_1\hat{y}_2}_{\hat{x}_2\hat{x}_1})^t\cdot \rmP,$$
where $(a_1\ot a_2)^t=a_1^t\ot a_2^t$ (here $a^t$ is the transpose of a matrix $a$) and $\rmP$ is given by \eqref{P-def},
and in the non-formal case the arguments $\hat{x}_i$ (resp., $\hat{y}_i$) take values in $\YY$ (resp., $\XX$).  
Then 

\noindent
(i) $\hat{r}^{\hat{x}_1\hat{x}_2}_{\hat{y}_1\hat{y}_2}$ is again a solution of the AYBE satisfying the skew-symmetry condition.
Furthermore, one has $\hat{\hat{r}}=r$.
%$r\mapsto \hat{r}$ is an involution.
If $\hat{x}_i$ are $\hat{y}_i$ are formal variables, and $r$ has an expansion of the form \eqref{r-expansion-eq} then 
so does $-\hat{r}^{\hat{x}_1\hat{x}_2}_{-\hat{y}_1,-\hat{y}_2}$.

\noindent
(ii) In the context of Theorem A, assume that $r^{x_1x_2}_{y_1y_2}$ corresponds to an $A_\infty$-structure on $\AA\ot R$.
Formally setting $\hat{X}=Y$ and $\hat{Y}=\hat{X}[1]$, $\hat{\th}_\alpha=\eta_\alpha$, $\hat{\eta}_\beta=\th_\beta$, etc.,
we get a new $A_\infty$-structure on $\AA\ot R$. Then the corresponding formal solution of the AYBE is precisely $\hat{r}$.
\end{prop}

\Pf . (i) Note that the map
$$\mathrm{Mat}_n(\k) \otimes \mathrm{Mat}_n(\k) \to \mathrm{Mat}_n(\k) \otimes
\mathrm{Mat}_n(\k) : x\mapsto \rmP\cdot x\cdot \rmP$$
is just the permutation of factors. Thus, the skew-symmetry condition can be rewritten as
$$\rmP r^{x_1x_2}_{y_1y_2}\rmP=-r^{x_2x_1}_{y_2y_1}.$$
Passing to the transposed matrices and making the substitution $x_i=\hat{y}_i$, $y_2=\hat{x}_1$, $y_1=\hat{x}_2$, 
we derive
that $\hat{r}$ also satisfies \eqref{skew-sym-eq}. The fact that $\hat{\hat{r}}=r$ follows from the identity $\rmP\cdot \rmP=1\ot 1$.

The AYBE equation for $\hat{r}$ is the following equation for $r$:
$$(r^{x_2x_3}_{y_1y_2})^{12,t}\rmP^{12}(r^{x_1x_3}_{y_3y_1})^{13,t}\rmP^{13}-
(r^{x_1x_2}_{y_3y_1})^{23,t}\rmP^{23}(r^{x_2x_3}_{y_3y_2})^{12,t}\rmP^{12}+
(r^{x_1x_3}_{y_3y_2})^{13,t}\rmP^{13}(r^{x_1x_2}_{y_2y_1})^{23,t}\rmP^{23}=0,$$
where we set $x_i=\hat{y}_i$, $y_i=\hat{x}_i$.
%$$(r^{12}(v,-u'))^tP^{12}r^{13}(v+v',u+u')^tP^{13}-
%r^{23}(v',u+u')^tP^{23}r^{12}(v,u)^tP^{12}+r^{13}(v+v',u)^tP^{13}r^{23}(v',u')^tP^{23}=0.$$
Using the skew-symmetry of $r$ (and the fact that $\rmP^t=\rmP$) we can rewrite this as
$$-\rmP^{12}(r^{x_3x_2}_{y_2y_1})^{12,t}(r^{x_1x_3}_{y_3y_1})^{13,t}\rmP^{13}+
\rmP^{23}(r^{x_2x_1}_{y_1y_3})^{23,t}(r^{x_2x_3}_{y_3y_2})^{12,t}\rmP^{12}-
\rmP^{13}(r^{x_3x_1}_{y_2y_3})^{13,t}(r^{x_1x_2}_{y_2y_1})^{23,t}\rmP^{23}=0.$$
%$$-P^{12}r^{12}(-v,u')^tr^{13}(v+v',u+u')^tP^{13}+
%P^{23}r^{23}(-v',-u-u')^tr^{12}(v,u)^tP^{12}-
%P^{13}r^{13}(-v-v',-u)^tr^{23}(v',u')^tP^{23}=0,$$
or passing to the transpose,
$$-\rmP^{13}(r^{x_1x_3}_{y_3y_1})^{13}(r^{x_3x_2}_{y_2y_1})^{12}\rmP^{12}+
\rmP^{12}(r^{x_2x_3}_{y_3y_2})^{12}(r^{x_2x_1}_{y_1y_3})^{23}\rmP^{23}-
\rmP^{23}(r^{x_1x_2}_{y_2y_1})^{23}(r^{x_3x_1}_{y_2y_3})^{13}\rmP^{13}=0.$$
%$$-P^{13}r^{13}(v+v',u+u')r^{12}(-v,u')P^{12}+
%P^{12}r^{12}(v,u)r^{23}(-v',-u-u')P^{23}-
%P^{23}r^{23}(v',u')r^{13}(-v-v',-u)P^{13}=0.$$
Using the fact that $\rmP\cdot \rmP=1$, we can rewrite this as
$$-(r^{x_1x_3}_{y_3y_1})^{13}(r^{x_3x_2}_{y_2y_1})^{12}+
\rmP^{13}\rmP^{12}(r^{x_2x_3}_{y_3y_2})^{12}(r^{x_2x_1}_{y_1y_3})^{23}\rmP^{23}\rmP^{12}-
\rmP^{13}\rmP^{23}(r^{x_1x_2}_{y_2y_1})^{23}(r^{x_3x_1}_{y_2y_3})^{13}\rmP^{13}\rmP^{12}=0.$$
%$$-r^{13}(v+v',u+u')r^{12}(-v,u')+P^{13}P^{12}r^{12}(v,u)r^{23}(-v',-u-u')P^{23}P^{12}-
%P^{13}P^{23}r^{23}(v',u')r^{13}(-v-v',-u)P^{13}P^{12}=0.$$
Now we use the identities $\rmP^{13}\rmP^{12}=\rmP^{12}\rmP^{23}$, $\rmP^{13}\rmP^{23}=\rmP^{12}\rmP^{13}$ and the fact that
$x\mapsto \rmP^{ij}x\rmP^{ij}$ acts as a transposition $(ij)$, to rewrite this as
$$-(r^{x_1x_3}_{y_3y_1})^{13}(r^{x_3x_2}_{y_2y_1})^{12}+(r^{x_2x_3}_{y_3y_2})^{23}(r^{x_2x_1}_{y_1y_3})^{31}
-(r^{x_1x_2}_{y_2y_1})^{12}(r^{x_3x_1}_{y_2y_3})^{32}=0.$$
%$$-r^{13}(v+v',u+u')r^{12}(-v,u')+r^{23}(v,u)r^{31}(-v',-u-u')-r^{12}(v',u')r^{32}(-v-v',-u)=0.$$
Swapping $2$ and $3$ we get the equation
$$-(r^{x_1x_3}_{y_3y_1})^{12}(r^{x_3x_2}_{y_2y_1})^{13}+(r^{x_2x_3}_{y_3y_2})^{32}(r^{x_2x_1}_{y_1y_3})^{21}
-(r^{x_1x_2}_{y_2y_1})^{13}(r^{x_3x_1}_{y_2y_3})^{23}=0.$$
Using the skew-symmetry to rewrite the middle summand we get the equation obtained from the AYBE by the change of variables
$$(x_1,x_2,x_3;y_1,y_2,y_3)\mapsto (x_3,x_1,x_2;y_2,y_3,y_1).$$
%$$(u,u',v,v')\mapsto (v',-v-v',u+u',-u).$$

\noindent
(ii) This follows from the cyclic symmetry. Namely, if we set $y_1=\hat{x}_2$, $y_2=\hat{x}_1$, $x_i=\hat{y}_i$, for $i=1,2$,
and use the corresponding identification between the twisted objects, we see that the new solution comes from
the Massey products for the composable arrows $Y_2\to X_1\to Y_1\to X_2$. Furthermore, it is given by
$$\sum_{\a,\a',\b,\b'}\lan \m_3^t(\eta^{12}_{\b},\th^{11}_\a,\eta^{21}_{\b'}),\th^{22}_{\a'}\ran\cdot e_{\a'\b}\ot e_{\a\b'}+\ldots,$$
where the other terms are standard singular parts.
Using the cyclic symmetry we can rewrite this as
$$\sum_{\a,\a',\b,\b'}\lan \m_3^t(\th^{22}_{\a'},\eta^{12}_{\b},\th^{11}_\a),\eta^{21}_{\b'}\ran\cdot e_{\a'\b}\ot e_{\a\b'}+\ldots.$$
This matches the formula for $(r^{x_1x_2}_{y_1y_2})^t\cdot \rmP$ due to the identity
$$(e_{\b'\a'}\ot e_{\b\a})^t\cdot \rmP=e_{\a'\b}\ot e_{\a\b'}.$$
\ed

\begin{cor}
    Let $r(u,v)$ be a solution of the AYBE with values in $\mathrm{Mat}_n(\k) \otimes
    \mathrm{Mat}_n(\k)$, satisfying the skew-symmetry condition \eqref{AYBE-skew-eq}. Then $r(v,u)^t\cdot \rmP$
is again a solution of the AYBE satisfying the skew-symmetry condition.
\end{cor}

\Pf . Apply Proposition \ref{propsymm}(i) to $r^{x_1x_2}_{y_1y_2}=r(x_1-x_2,y_1-y_2)$.
  \ed
  
Recall that the nondegeneracy condition on solutions of the AYBE imposed in \cite{pol09}
is that the tensor $r(u,v)\in \mathrm{Mat}_n(\k) \ot \mathrm{Mat}_n(\k)$ is nondegenerate for generic $(u,v)$.
Now we are going to use the above involution to show that the pole conditions for $r(u,v)$, imposed in the classification
result of \cite{pol09}, are implied by the following stronger nondegeneracy condition, involving $r(u,v)$ and $r(u,v)^t\cdot \rmP$.

\begin{defi} \label{strdegdef} Let us say that an
    $\mathrm{Mat}_n(\k)\otimes \mathrm{Mat}_n(\k)$-valued function $r(u,v)$, meromorphic in a neighborhood of $(0,0)$, 
is {\it strongly nondegenerate} if the tensors $r(u,v)$ and $r(u,v)^t\cdot \rmP$ are nondegenerate for generic $(u,v)$.
\end{defi}

\begin{prop} \label{strdeg}  Assume $N>1$. Let $r(u,v)$ be a strongly nondegenerate skew-symmetric 
solution of the AYBE (meromorphic in a neighborhood of $(0,0)$).
Then $r(u,v)$ has a simple pole at $u=0$ (resp., $v=0$) with the polar term
$c\cdot \frac{1\ot 1}{u}$ (resp., $c'\cdot \frac{\rmP}{v}$), where $c$ and $c'$ are nonzero constants.
\end{prop}

\Pf . First, we claim that the involution $r(u,v)\mapsto \hat{r}(u,v):=r(v,u)^t\cdot \rmP$ 
on skew-symmetric solutions of the AYBE preserves the notion of strong
nondegeneracy. Indeed, this immediately follows from the observation that $\rmP\cdot r\cdot \rmP=r^{21}$, so it
is nondegenerate if and only if $r$ is nondegenerate.

Next, given $r(u,v)$, a strongly nondegenerate skew-symmetric solution of the AYBE, we claim that $r(u,v)$
necessarily has a pole at $u=0$. Indeed, assume $r(u,v)$ is regular at $u=0$. Then by \cite[Lem.\ 1.2]{pol09},
we have an expansion $r(u,v)=r_0(v)+ur_1(v)+\ldots$, where $r_0(v)=r(u,0)$ is a nondegenerate skew-symmetric 
  solution of the 
AYBE. Then by \cite[Thm.\ 0.2]{pol09}, $r_0(v)$ has a pole at $v=0$, hence, $r(u,v)$ has a pole at $v=0$. By 
\cite[Lem.\ 1.3]{pol09}, this implies that $r(u,v)$ has a simple pole at $v=0$ with the polar part $c\cdot\frac{\rmP}{v}$.
Hence, $\hat{r}(u,v)$ has a simple pole at $u=0$ with the polar part $c\cdot \frac{1\ot 1}{u}$.
Since $\hat{r}(u,v)$ is still
a nondegenerate skew-symmetric solution of the AYBE, by \cite[Lem.\ 1.5]{pol09}, $\hat{r}(u,v)$ has a simple pole at $v=0$.
Equivalently, $r(u,v)$ has a simple pole at $u=0$, which is a contradiction.

Thus, we know that $r(u,v)$ has a pole at $u=0$, or equivalently, $\hat{r}(u,v)$ has a pole at $v=0$. By \cite[Lem.\ 1.3]{pol09},
this implies that $\hat{r}(u,v)$ has a simple pole at $v=0$ with the polar part $c\cdot\frac{\rmP}{v}$. Hence,
$r(u,v)$ has a simple at $u=0$ with the polar part $c\cdot \frac{1\ot 1}{u}$. Now the assertion follows from \cite[Lem.\ 1.5]{pol09}.
\ed

\subsection{From algebraic/analytic to formal solutions of the general AYBE}\label{alg-form-AYBE-sec}

In this section we want to consider the solutions of the general AYBE arising, as described in Introduction,
from two algebraic families of objects $\XX$ and $\YY$. We want to show how to pass from these solutions
to the corresponding formal solutions associated to picking one object in each family.

We will use the formalism from \cite[Sec.\ 1]{Seidel-flux}
concerning families of objects in $A_\infty$-categories.

Let $\mathscr{A}$ be an $A_\infty$-category over $\k$, and let $\XX$ and $\YY$ be smooth affine curves over
$\k$, such that we have
perfect families of $\mathscr{A}$-modules $\MM$ and $\NN$ parametrized by $\XX$ and $\YY$.
We assume that for $x\neq x'$ (resp., $y\neq y'$) one has $\Hom^*(M_x,M_{x'})=0$ (resp., $\Hom^*(N_y,N_{y'})=0$),
that each $M_x$ (resp., $N_y$) is $1$-spherical, and that $\Hom(M_x,N_y)$ are concentrated in degree $0$.
Furthermore, we assume that $\sHom^*(p_{\XX}^*\MM,p^*_{\YY}\NN)$ is a vector bundle over $\XX\times\YY$.

Recall (see \cite[(1h)]{Seidel-flux}) that one can associate with the families $\MM$ and $\NN$ the deformation classes
$$Def(\MM)\in \Om^1_\XX\ot \Hom^1(\MM,\MM), \ \ 
Def(\NN)\in \Om^1_\YY\ot \Hom^1(\NN,\NN).$$

Let $\UU\sub \XX^2\times \YY^2$ be the complement to the diagonals 
$\De_\XX\times\YY^2\cup \XX^2\times\De_\YY$. Then over $\UU$ we have the induced families 
$\MM(x_1)$, $\MM(x_2)$, $\NN(y_1)$ and $\NN(y_2)$ (pull-backs from the families over $\XX$ and $\YY$), such that 
$$\sHom^*(\MM(x_1),\MM(x_2))=\sHom^*(\NN(y_1),\NN(y_2))=0.$$
Thus, we have a well defined triple Massey product 
\begin{align*}
    &\mathrm{MP}(x_1,x_2;y_1,y_2):\sHom^0(\MM(x_1),\NN(y_1))\ot\sHom^1(\NN(y_1),\MM(x_2))\ot\sHom^0(\MM(x_2),\NN(y_1))\to \\
&\sHom^0(\MM(x_1),\NN(y_2)),
\end{align*}
which is a morphism of vector bundles over $\UU$.

On the other hand, let us fix points $x_0\in\XX$ and $y_0\in\YY$, and let us fix generators
$\xi_{x_0}\in\Hom^1(M_{x_0},M_{x_0})$ and $\xi_{y_0}\in \Hom^1(N_{y_0},N_{y_0})$. Assume that
$$Def(\MM)|_{x_0}\in \Om^1_\XX|_{x_0}\ot\Hom^1(M_{x_0},M_{x_0})\simeq \Om^1_\XX|_{x_0}$$
is nonzero, and similarly $Def(\NN)|_{y_0}\in \Om^1_\YY|_{y_0}$ is nonzero.
Then we can choose a formal parameter $t$ near $x_0$ on $\XX$ (resp., parameter $s$ near $y_0$ on $\YY$)
such that $Def(\MM)=dt$ in a formal neighborhood of $x_0$ (resp., $Def(\NN)=ds$ near $y_0$).

Let $\hat{\MM}_{x_0}$ be the family over $k[[t]]$ obtained from $\MM$ by restriction to a formal disk around $x_0$,
and let $\hat{\NN}_{y_0}$ be the similar family over $k[[s]]$.
On the other hand, we have twisted objects $(M_{x_0},t\xi_{x_0})$ and $(N_{y_0},s\xi_{y_0})$, over $k[[t]]$ and $k[[s]]$, respectively. These twisted objects produce the same deformation classes $dt$ and $ds$, so
by the proof of \cite[Prop.\ 1.21]{Seidel-flux}, we derive the existence of quasi-isomorphisms of families
$$\hat{\MM}_{x_0}\simeq (M_{x_0},t\xi_{x_0}), \ \ \hat{\NN}_{y_0}\simeq (N_{y_0},s\xi_{y_0}).$$

By the functoriality of Massey products, this implies that the formal expansion of the Massey
product $\mathrm{MP}(x_1,x_2,y_1,y_2)$
near $x_0$ and $y_0$, is equal to the triple Massey product considered in the proof of Theorem A.

Note that $Def(\MM)|_{x_0}$ can be identified with the usual class of the first-order deformation of $M_{x_0}$,
associated with $\MM$. In particular, in the situation when $\MM$ is a universal deformation of $M_{x_0}$ then
$Def(\MM)|_{x_0}$ is nonzero. This is the situation that occurs when we consider families of simple vector bundles
(or structure sheaves of points) on elliptic curves and their degenerations, as in \cite{pol02}, \cite{pol09}.
In the case of families of Lagrangians in Fukaya category, we have a similar picture, with algebraic
families replaced by analytic families.

\section{Trigonometric solutions of the AYBE from symplectic geometry}
\label{symplectic}

\subsection{A square-tiled surface from Belavin-Drinfeld structures}

In this section, starting from an associative Belavin-Drinfeld structure $(S,C_1,C_2,A)$, we
construct a punctured Riemann surface $\Sigma$ together with a non-vanishing holomorphic one-form $\alpha \in \Gamma(C,\Omega^{1,0}_C)$. 

Let us begin with a finite set $S$ of $n$ elements, and two permutations $C_1, C_2 \in Aut(S) \cong \mathfrak{S}_n$ such that the subgroup $\langle C_1, C_2 \rangle \subset Aut(S)$ is transitive. Let $\mathbb{T}$ be the square torus $\C  / (\Z \oplus i \Z)$ and let $\mathbb{T}_0 = \mathbb{T} \setminus \{0 \}$ be the punctured torus. Let us also consider the oriented curves, $l_1,l_2: [0,1] \to \mathbb{T}$ defined by:
\[ l_1(t) = \frac{1-it}{2}   \ , \ \ \  l_2 (t) = \frac{t+i}{2}   \]
Let $p_0\in \T_0$ be the point $(1+i)/2$, which is the unique point in $l_1\cap l_2$.
Consider the $n$-fold (unramified) covering:
\[ \pi : \Sigma_0 \to \mathbb{T}_0 \]
corresponding to the subgroup $H \subset \pi_1(\mathbb{T}_0 , l_1 \cap l_2) = \langle l_1, l_2 \rangle \cong F_2$ given as the preimage of the stabilizer of a point $s_0\in S$
under the homomorphism $\rho : \pi_1(\mathbb{T}_0) \to Aut(S)$ defined by \[ \rho(l_1) = C_1 \ , \ \ \ \rho(l_2) = C_2. \]
Thus, we identify the set $S$ with the fiber $\pi^{-1}(p_0)$, so that the action of the generators $l_1$ and $l_2$
of the fundamental group $\pi_1(T_0,p_0)$ on the fibre $S$ is given by the permutations $C_1, C_2 : S
\to S$.
 We let $L_i = \pi^{-1}(l_i)$ be the multi-curves in $\Sigma_0$ covering the circles $l_i$, so that 
$L_1 \cap L_2 = \pi^{-1}(p_0)=S$. 

The assumption that the subgroup $\langle
C_1,C_2 \rangle \subset Aut(S)$ is transitive guarantees that $\Sigma_0$ is connected. The curves $L_1$ and $L_2$ are
connected if and only if both $C_1$ and $C_2$ are transitive
permutations. We will always require that $\Sigma_0$ is connected.
If $(S,C_1,C_2)$ comes from an associative Belavin-Drinfeld data, then $C_1, C_2$
are required to be transitive permutations, however this condition
is not strictly necessary in what follows. 

One can lift the flat metric on $\mathbb{T}_0$ to $\Sigma_0$. To visualise this
metric on $\Sigma_0$, let us now give a more geometric construction of the
covering map $\pi : \Sigma_0 \to \mathbb{T}_0$. Recall that $\mathbb{T}$ is
obtained from the unit square in $[0,1] \times [0,1] \subset \mathbb{R}^2$ by
identifying the opposite sides. $\mathbb{T}_0$ is obtained from this by
removing the corner point. Now let us take $n$ copies of the unit square (with
corners removed) labeled the set $S = \{1,\ldots, n\}$. Given
automorphisms $C_1, C_2 \subset \mathfrak{S}_n \cong Aut(S)$, construct a
surface $\Sigma_0$ as follows: 1) identify the right edge of the $i^{th}$
square with the left edge of $j^{th}$ square if $C_1(i)=j$; 2) identify the bottom edge of the $i^{th}$ square with the top edge of the $j^{th}$ square if
$C_2(i)=j$. It is because of this construction $\Sigma_0$ is called a \emph{square-tiled surface}. The
name was first suggested to Anton Zorich by Alex Eskin \cite{Zmiaikou}.

By the Riemann existence theorem (\cite[Sec. 4.2.2]{donaldson}), the surface $\Sigma_0$ can be
completed to a surface $\widehat{\Sigma}_0$ and the covering map extends to a branched covering map:
\[ \widehat{\pi}: \widehat{\Sigma}_0 \to \mathbb{T} \]
ramified along the origin $(0,0) \in \mathbb{T}$. The preimage of the origin, 
  \[ \{ p_1,p_2, \ldots p_b \} = \widehat{\pi}^{-1}(0) \]
consists of a number of points, which is equal to the number of cycles in the cycle decomposition of
the commutator $[C_1,C_2]$ into disjoint union of cycles of varying lengths (from $1$ to $n$).
Indeed, the curves $L_1$ and $L_2$ divide the surface $\widehat{\Sigma}_0$ into polygons, such that the point
$p_i$ is contained in a $(2e(p_i)+2)$-gon, where $e(p_i)$ is the ramification index of the point $p_i$. 

We let $b_k$ denote the number of $k$-cycles, so that we have $n=\sum_{k=1}^n k b_k$. We record the
following elementary computation, which follows from the above explicit description of $\widehat{\Sigma}_0$ as a union of $b$-polygonal regions, or also by the Riemann-Hurwitz formula.
\begin{prop} The number of points in $\widehat{\Sigma_0} \setminus \Sigma_0$ is equal to $b= \sum_{k=1}^n b_k$. The Euler characteristic of $\Sigma_0$ is $\chi(\Sigma_0)=-n$. Consequently, the genus $g$ is determined by the formula \[ \chi(\Sigma_0) = 2-2g -b = -n. \] 
In particular, $g=1$ if and only if $C_1$ and $C_2$ commute. \ed \end{prop} 

Finally, we need to incorporate the proper subset $A \subset S$ that appears in
an associative Belavin-Drinfeld structure $(S,C_1, C_2, A)$. This data enters
in determining a partial compactification of $\Sigma_0$. 

Namely, recall that $A$ is, by definition, a subset of the set of fixed points of the action of the
commutator $[C_1,C_2]$ on $S$. In terms of the branched covering map $\widehat{\pi}:
\widehat{\Sigma}_0 \to \mathbb{T}$, the set of fixed points of the commutator $[C_1,C_2]$ can be identified with the set of points $p_i$ in
the preimage $\widehat{\pi}^{-1}(0)$ which have ramification index $e(p_i)=1$. Using this we
can identify $A$ with a subset of points $p_i$ where the map $\widehat{\pi}$ is unramified. To be precise, an element $a \in A$ gives a square with corners \[ \{ a, C_2(a), C_1C_2(a),
C_2^{-1}C_1C_2(a) \}, \] which determines a point $p_{a} \in \widehat{\pi}^{-1}(0) $ of ramification
index 1 contained in this square. We define $\Sigma= \Sigma_A $ to be
the partial compactification $\Sigma_0 \cup \{ p_{a} | a\in A \}$. Note that the covering map
$\pi : \Sigma_0 \to \mathbb{T}_0$ extends to a local diffeomorphism:
\[ \pi : \Sigma \to \mathbb{T} \]
Hence, the  flat metric on $\mathbb{T}$ lifts to a flat metric on $\Sigma$ so as to make $\pi$ into
a local isometry. From now on, we will consider the square-tiled surface $\Sigma$ equipped with this
flat metric. Note that, for convenience, we always normalize the metric on $\mathbb{T}$ so that the length of the curves $l_1$
and $l_2$ are 1. 

Equivalently, we write $C$ for the unique Riemann surface structure on $\Sigma$ making $\pi: \Sigma \to \T$ into a holomorphic map. In this case, we equip $C$ with the one-form $\alpha = \pi^* dz$, the pullback of the standard non-vanishing holomorphic one-form on $\T$. 

\begin{ex}
    \label{ex:sqtiled}
Figure \ref{fig1} shows an example of this construction corresponding to $S = \{1,2,3,4\}, C_1
=(1,4,2,3), C_2 = (1,2,3,4)$. The red curve $L_1 \subset \Sigma_0$ and the blue curve $L_2 \subset \Sigma_0$ depict the preimages of the curves $l_1 , l_2
\subset \mathbb{T}_0$. The flat metric can be extended over the black labelled point without any
singularities. Thus, we can choose $A$ to be either empty or include the unique black
labelled point, which corresponds to $\{3\}$ - the unique fixed point of $[C_1,C_2]$. If $A=\{3 \}$,
then correspondingly, we compactify $\Sigma_0$ by filling in the puncture labelled black.

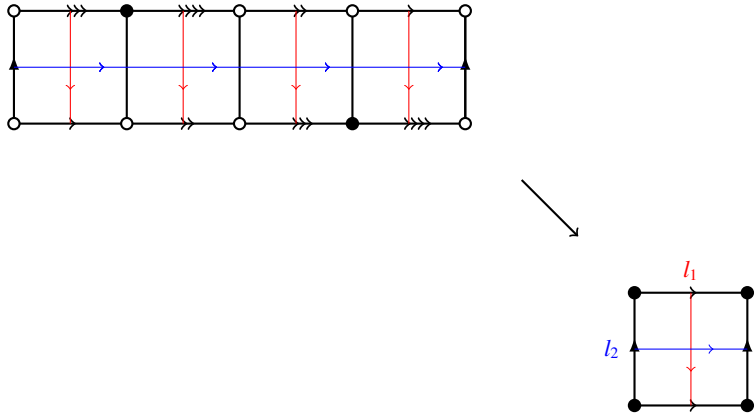
\begin{figure}[htb!]
\centering
\begin{tikzpicture} [scale=1.5]

        \tikzset{->-/.style={decoration={ markings,
	    mark=at position #1 with {\arrow{>}}},postaction={decorate}}}
        \tikzset{->>-/.style={decoration={ markings,
	    mark=at position #1 with {\arrow{>>}}},postaction={decorate}}}
        \tikzset{->>>-/.style={decoration={ markings,
	    mark=at position #1 with {\arrow{>>>}}},postaction={decorate}}}
       \tikzset{->>>>-/.style={decoration={ markings,
	    mark=at position #1 with {\arrow{>>>>}}},postaction={decorate}}}
        \tikzset{-|>-/.style={decoration={ markings,
	    mark=at position #1 with {\arrow{latex}}},postaction={decorate}}}

	\draw[thick, ->-=.55] (-15,0) -- (-14,0);
    \draw[thick, ->>-=.6] (-14,0) -- (-13,0);
    \draw[thick, ->>>-=.65] (-13,0) -- (-12,0);
        \draw[thick, ->>>>-=.7] (-12,0) -- (-11,0);
        
        \draw[thick, ->>>-=.65] (-15,1) -- (-14,1);
        \draw[thick, ->>>>-=.7] (-14,1) -- (-13,1);
        \draw[thick, ->>-=.6] (-13,1) -- (-12,1);
        \draw[thick, ->-=.55] (-12,1) -- (-11,1);
       
	\draw[thick, -|>-=.6] (-15,0) -- (-15,1);
	\draw[thick, -|>-=.6] (-11,0) -- (-11,1);

	\draw[thick] (-14,0) -- (-14,1);
	\draw[thick] (-13,0) -- (-13,1); 
	\draw[thick] (-12,0) -- (-12,1); 
	\draw[thick] (-11,0) -- (-11,1); 
         
	\draw[red, ->-=.7] (-14.5,1) -- (-14.5,0); 
        \draw[red, ->-=.7] (-13.5,1) -- (-13.5,0); 
	\draw[red, ->-=.7] (-12.5,1) -- (-12.5,0); 
        \draw[red, ->-=.7] (-11.5,1) -- (-11.5,0); 

        \draw[blue, ->-=.2, ->-=.45, ->-=.7, ->-=.95] (-15,0.5) -- (-11,0.5); 
	\draw[thick,->] (-10.5,-0.5) -- (-10,-1);

	\draw [thick, ->-=.55] (-9.5,-1.5) -- (-8.5, -1.5);
        \draw [thick, ->-=.55] (-9.5,-2.5) -- (-8.5,-2.5);
        \draw [thick, -|>-=.6] (-9.5,-2.5) -- (-9.5, -1.5);
        \draw [thick, -|>-=.6] (-8.5,-2.5) -- (-8.5,-1.5);
       
	\draw[red, ->-=.7] (-9,-1.5) -- (-9,-2.5); 
        \draw[blue, ->-=.7] (-9.5,-2) -- (-8.5,-2); 

        \node[blue] at (-9.7,-2)   {\footnotesize $l_2$};
        \node[red] at (-9,-1.3)   {\footnotesize $l_1$};

        \draw[thick, fill=white] (-15,0) circle(.05); 
        \draw[thick, fill=white] (-14,0) circle(.05); 
        \draw[thick, fill=white] (-13,0) circle(.05); 
        \draw[thick, fill=black] (-12,0) circle(.05); 
        \draw[thick, fill=white] (-11,0) circle(.05); 
        \draw[thick, fill=white] (-15,1) circle(.05); 
        \draw[thick, fill=black] (-14,1) circle(.05); 
        \draw[thick, fill=white] (-13,1) circle(.05); 
        \draw[thick, fill=white] (-12,1) circle(.05); 
        \draw[thick, fill=white] (-11,1) circle(.05); 
    
        \draw[thick, fill=black] (-9.5,-1.5) circle(.05); 
        \draw[thick, fill=black] (-9.5,-2.5) circle(.05); 
        \draw[thick, fill=black] (-8.5,-1.5) circle(.05); 
        \draw[thick, fill=black] (-8.5,-2.5) circle(.05); 
    
\end{tikzpicture}
	\caption{A square-tiled surface}
	\label{fig1}
\end{figure}

\end{ex} 
\begin{figure}[htb!]
\centering
\begin{tikzpicture}
\draw (0,0) node[inner sep=0] {\includegraphics[scale=0.7]{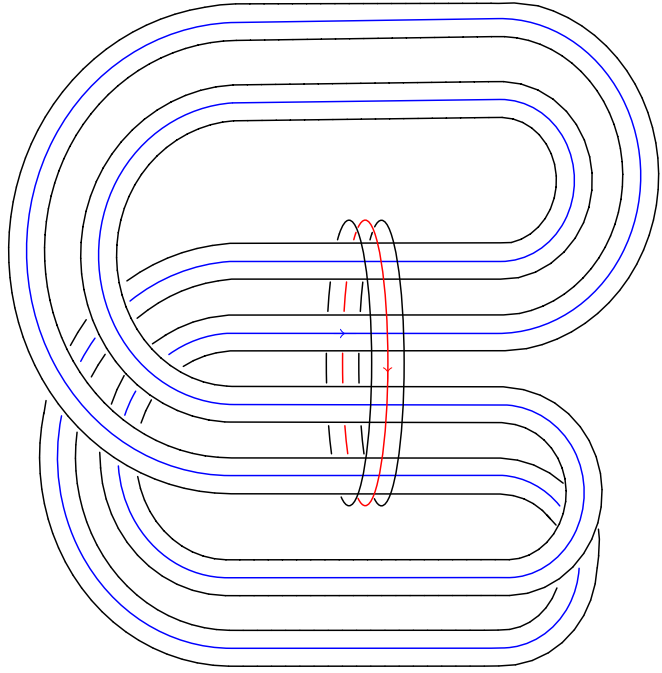}};
\tikzset{->-/.style={decoration={ markings,
	    mark=at position #1 with {\arrow{>}}},postaction={decorate}}}
\draw[red, ->-=.1] (0.717,-0.5) -- (0.717,-0.65); 
\draw[blue, ->-=.5] (0,0) -- (0.3,0); 
\end{tikzpicture}
\caption{$S = \{1,2,3,4\}, C_1 =(1,2,3,4), C_2 = (1,3,2,4)$}
\label{fig2}
\end{figure}
\begin{rem} Topologically the punctured torus $\mathbb{T}_0$ can be seen as the plumbing of two
    cotangent bundles $T^*l_i \cong T^*S^1$ at one point $l_1 \cap l_2$. Similarly, one can
    construct the surface $\Sigma_0$ as the plumbing of $T^*L_1$ and $T^*L_2$ at $n$-points
    corresponding to $L_1 \cap L_2$. The Figure \ref{fig2} illustrates this construction in the case
    of $S=\{1,2,3,4\}$, $C_1=(1,2,3,4)$,
 $C_2=(1,3,2,4)$. \end{rem}

\subsection{A Fukaya category from Belavin-Drinfeld structures}
Let $(C,\alpha)$ be the square-tiled surface obtained from an associative
Belavin-Drinfeld structure as above. Let $\Sigma$ be the topological surface
underlying $C$. The square $\Omega = \alpha \otimes \alpha \in \Gamma(C,(\Omega_C^{1,0})^{\otimes 2})$ determines a
non-vanishing quadratic form, which gives a flat Riemannian metric $|\Omega|$
on $\Sigma$ and a horizontal foliation of tangent vectors $v$ with $\Omega(v,v)
>0$. The Riemannian metric determines an area form\footnote{Note that our symplectic form
$\omega$ is not convex at infinity. This is usually required in setting up Floer theory in order to
ensure that a maximum principle holds which guarantees that pseudo-holomorphic disks remain in a
compact region. However,
in dimension 2, this holds for topological reasons. Alternatively, one could modify
$\omega$ near infinity to make it convex. Either way, the outcome is unchanged and we will simply use the area form $\omega$.} $\omega$ and the horizontal
foliation determines a grading structure on $\Sigma$, i.e a section of the
projectivized tangent bundle of $S$, which we view as an unoriented line field $l \subset T(\Sigma)$. We note that such line fields form a torsor for $C^{\infty}(\Sigma,
\mathbb{R}P^1)$, and the connected components of this group can be identified with $H^1(\Sigma; \mathbb{Z})$. 

To work over $\C$, one works with exact Fukaya categories as in
\cite{seidelbook}. Thus, we will need to choose a primitive $\theta$ for $\omega$, which exists since $\Sigma$ is non-compact. We choose this so that the Lagrangians $L_1$ and $L_2$ are
exact. One can arrange this as follows: Choose any primitive $\theta_0$ for $\omega$; find a closed
$1$-form $\sigma_0$ such that $\sigma_0 ([L_i]) = \int_{L_i} \theta_0$, which exists
since $L_i$ give independent non-trivial homology classes in $H_1(\Sigma)$; and let $\theta =
\theta_0 - \sigma_0$.  We also normalize the area form so that the geodesics $L_1$ and $L_2$ have length 1.

We can now form a $\mathbb{Z}$-graded triangulated $1$-Calabi-Yau $\C$-linear $A_\infty$
category, the Fukaya category $\mathcal{F}(\Sigma)$ of $(\Sigma, d\theta,
\Omega)$. The objects of $\mathcal{F}(\Sigma)$ are closed, exact, oriented curves
$L$ equipped with grading structures and a rank $1$ local system $\xi
\to L$. 
%whose fibre is $\C$ and holonomy is in $\C^\times$.
Recall that a grading structure on a curve $L$ means a choice of a homotopy
class of a path from $TL$ to the line field $l_{|L}$. If $x \in L \cap L'$ is a transverse
intersection point, the grading $|x|$ is given by $\lfloor \alpha/\pi
\rfloor+1$ where $\alpha$ is the net rotation from $T_x L \to l_x \to T_x
L'$. This lifts the $\mathbb{Z}_2$-grading on the intersection points given by
$\lfloor (L \cdot_x L') / \pi \rfloor +1$ where $L \cdot_x L'$ is the local
algebraic intersection number of $L$ and $L'$ at $x$ associated to
orientations of $L$ and $L'$. 

Note also that on a circle there are precisely two spin structures corresponding to connected and disconnected double coverings of the circle. We implicitly fix a spin structure on each closed, exact Lagrangian $L \subset \Sigma$. Changing the spin structure by the action of $H^1(S^1;\mathbb{Z}_2) = \mathbb{Z}_2$ is equivalent to modifying the monodromy of
the local system $\xi \to L$ by the action of $\{ \pm 1 \} \subset \C^\times$. Therefore, the effect of changing the choice of spin structure on $L$ can be achieved by modifying the $\C^\times$ local system $\xi$. Spin structures enter in
defining orientations of various moduli spaces of holomorphic curves and they play a role in
determining the signs in various counts.  In the case of Fukaya categories of 2-dimensional
surfaces, which is the only situation considered in this paper, there is a combinatorial method
given in \cite[Sec.  7]{seidelgenus2} that allows us to compute these signs. Throughout, in our
explicit computations, we follow this method to determine the signs without giving further
explanation. 

The morphism spaces in the Fukaya category are given by Floer cochain complexes:
\[ CF^*((L_1, \xi_1), (L_2, \xi_2)) = \bigoplus_{x \in L_1 \cap L_2}
\mathrm{hom}_{\C} (\xi_1|_{x}, \xi_2|_{x}) \]

For brevity, we often suppress the local systems $\xi_i$ from the notation. The $A_\infty$-structure
comprises a collection of maps: \[ \m_k : CF(L_{k-1},
	L_k) \otimes \ldots \otimes CF(L_0,L_1) \to CF(L_0, L_k)[2-k] \] 
For $p_i \in L_{i-1} \cap L_{i}$ and $p_0 \in L_0 \cap L_k$, 
%$\rho_i \in \mathrm{hom}_{\C}(\xi_{i-1}|_{p_i}, \xi_{i}|_{p_i})$  $\rho_0 \in \mathrm{hom}_{\C} (\xi_0|_{p_0}, \xi_{k}|_{p_0})$, the matrix coefficient $ \langle \rho_{k},\ldots,\rho_1 | \rho_0 \rangle$ of 
the components of these maps involving $\mathrm{hom}_{\C}(\xi_{i-1}|_{p_i}, \xi_{i}|_{p_i})$ are
	defined by counting holomorphic disks with $(k+1)$-boundary punctures
	such that the boundary components are mapped to $(L_0,L_1,\ldots L_k)$
	in the cyclic order. Let us denote the moduli space of such
	pseudoholomorphic disks $u$ in the homotopy class $[u]$ by
	$\mathcal{M}(p_k,p_{k-1},\ldots,p_1,p_0; [u])$. If the index of $[u]$ is fixed to be $2-k$ and
    the regularity is arranged then Gromov-Floer compactness ensures that this moduli space is a
    finite set of points, which we can then count (with signs). 
    For $\rho_i \in \mathrm{hom}_{\C}(\xi_{i-1}|_{p_i}, \xi_{i}|_{p_i})$,
we set
 \begin{equation}\label{matrixcoef} \m_k(\rho_k,\ldots,\rho_1) = \sum_{\substack{ [u]: \text{ind}([u])=2-k}  } \# \mathcal{M}(p_k,p_{k-1},\ldots,p_1,p_0; [u])\cdot \text{hol}_{\partial u}  \in  \mathrm{hom}_{\C} (\xi_0|_{p_0}, \xi_{k}|_{p_0}),\end{equation} 
where the term $\text{hol}_{\partial{u}}$ is defined as follows. The boundary component of $u$ mapping to
$L_i$ gives isomorphisms $\xi_i|_{p_i} \to \xi_i|_{p_{i+1}}$. Therefore, given elements $\rho_i \in
\mathrm{hom}_\C(\xi_{i-1}|_{p_i}, \xi_i|_{p_i})$, using the isomorphisms provided by $u$, we can
construct the composition: \[ \text{hol}_{\partial{u}} = \rho_k \circ \ldots \circ \rho_1 \in
\mathrm{hom}_{\C}(\xi_0|_{p_0}, \xi_k|_{p_0}).  \]
Note that exactness ensures that the sum in Equation (\ref{matrixcoef}) is finite, hence is well-defined.     

In practice, when we do explicit computations, we will mark points by $\star$ on each Lagrangian
circle, and the contribution of a holomorphic disk will be weighted by the holonomy factor each time
the boundary of the disk passes through the marked point. 

We also note that it was proven by Fukaya \cite{fukayacyclic} that the Fukaya category of compact Lagrangians over $\mathbb{R}$ (equivalently, over any field $\k$ of characteristic 0) has a model with a strictly cyclic $A_\infty$-structure. The existence of such a cyclic structure is important for the applicability of Theorem A to $A_\infty$ algebras that we compute from Fukaya categories below. On the other hand, for the purpose of computation of triple Massey products, we can use any model of the Fukaya category as triple Massey products are homotopy invariant notions. We find it convenient to use the model of the Fukaya category as given in \cite[Sec. 7]{seidelgenus2}. 

\subsection{Constructing solutions to the AYBE via Massey products in $\mathcal{F}(\Sigma)$} 
\label{Fukaya}

As was shown in \cite{pol09}, with every associative Belavin-Drinfeld structure $(S,C_1,C_2,A)$ one
can associate a trigonometric solution of the AYBE. A slightly different looking
definition of an associative Belavin-Drinfeld structure was used in \cite{pol09}. In the next lemma
we show the equivalence of the Definition \eqref{BD} with the definition of the associative
Belavin-Drinfeld structure in \cite{pol09}.

Let $S$ be a finite set of $n$ elements. We denote a transitive permutation as a map $C: S \to S$ and we write $\Gamma_C := \{ (s, C(s)) | s
\in S \} \subset S \times S$ for its graph. 

\begin{lem} Let $S$ be a set equipped with a pair of transitive permutations $C_1, C_2 : S \to S$.
Then to give a proper subset $A \subset S$, such that $(S,C_1,C_2,A)$ is
an associative Belavin-Drinfeld structure, is
equivalent to giving
a pair of proper subsets $\Gamma_1, \Gamma_2 \subset \Gamma_{C_1}$ such that $(C_2 \times C_2) \Gamma_1 = \Gamma_2$.
\end{lem}

\Pf . We set 
\begin{align*} 
	\Gamma_1 = \{ (a, C_1(a)) | a \in A \}, \ \ \ \Gamma_2 = \{ (C_2(a), C_1(C_2(a)) | a \in A. \} 
\end{align*} 
One immediately sees that the condition
\[  (C_2 \times C_2) (\Gamma_1) = \Gamma_2 \] 
is equivalent to $C_1C_2(a)=C_2C_1(a)$ for every $a\in A$.
\ed

We prefer the form given in Definition \ref{BD} as it makes the symmetry with respect to switching
$C_1$ and $C_2$ more clear (cf. Proposition \ref{propsymm}). In examples, it may be convenient to identify $S = \{1, \ldots, n \}$ such that $C_1(i) = i+1$ (modulo $n$). One then thinks of $C_2$ as an $n$-cycle in the symmetric group $\mathfrak{S}_n$.

The commutator $[C_1, C_2] \in \mathfrak{A}_n \subset \mathfrak{S}_n$ plays a
special role in the definition as the elements of the set $A$ correspond to a
subset of the set of fixed points of the commutator $[C_1, C_2]$. We remark
that it can be proven by induction that any element of the alternating group
$\mathfrak{A}_n$ arises as the commutator $[C_1,C_2]$ of two $n$-cycles in
$\mathfrak{S}_n$ (see Prop. 4 \cite{huse}). 

Let $A_S$ denote
the algebra of endomorphisms of the $\C$-vector space with the basis $(\be_i)_{i\in S}$,
so that $A_S\simeq\Mat_n(\C)$,  where $n=|S|$. We denote by $e_{ij}\in A_S$ the endomorphism
defined by $e_{ij}(\be_k)=\de_{jk}\be_i$.
The solution of the AYBE associated with $(S,C_1,C_2,A)$ is given by
\begin{equation}
\begin{aligned}
	\label{trig-AYBE-eq} 
    &r(u,v)  = \frac{1}{\exp(u)-1}\sum_{i} e_{ii}\ot e_{ii} + \frac{1}{1-\exp(-v)}\sum_i e_{ii}\ot e_{ii}  \\        
       & +      \frac{1}{\exp(u)-1}\sum_{0< k<n, i}\exp(\frac{ku}{n})e_{C_1^k(i),C_1^k(i)}\ot e_{ii}
       + \frac{1}{\exp(v)-1}\sum_{0<m<n, i}\exp(\frac{mv}{n})e_{i,C_2^m(i)}\ot e_{C_2^m(i),i}  \\
       & +\sum_{0<k,0<m; a\in A(k,m)}
\Big\{ \exp(-\frac{ku+mv}{n})e_{C_2^m(a),a}\ot
e_{C_1^k(a),C_1^kC_2^m(a)}-\exp(\frac{ku+mv}{n})e_{C_1^k(a),C_1^kC_2^m(a)}\ot e_{C_2^m(a),a}
\Big\},  \end{aligned}\end{equation}
    where we denote by $A(k,m)\sub A$ the set of all $a\in A$ such that
$C_1^iC_2^j(a)\in A$ for all $0\le i< k, 0\le j< m$. One can easily check that for $a\in A(k,m)$ one has
$C_1^kC_2^m(a)=C_2^mC_1^k(a)$. Note also that $A(k,m)$ can be nonempty only if $k<n$ and $m<n$ (since $A$
is a proper subset of $S$), so our formula is equivalent to that of \cite[Thm.\ 0.1]{pol09}.

Let us denote by $\pr:\Mat_n(\C)\to\ssl_n(\C)$ the projection along $\C\cdot 1$.
Let $r(u,v)$ be a unitary solution of the AYBE such that the Laurent expansion of $r$ at $u=0$ has
form 
\begin{equation}\label{Laurent-eq}
r(u,v)=\frac{1\ot 1}{u}+r_0(v)+ur_1(v)+\ldots.
\end{equation}
Then one can show that $(\pr\ot\pr)r_0(v)$ is a unitary solution of the CYBE, nondegenerate if $r(u,v)$ was nondegenerate.

One of the main results of \cite{pol09} is that 
every nondegenerate unitary solution of the AYBE for $A=\Mat_n(\C)$ (where $n>1$), 
such that the Laurent expansion of $r$ at $u=0$ has
form \eqref{Laurent-eq} and $(\pr\ot\pr)r_0(v)$ is a trigonometric solution of the CYBE,  is equivalent to one of the solutions \eqref{trig-AYBE-eq}.

We will next show that the above solutions to the AYBE can be recovered from Massey products in $\mathcal{F}(\Sigma)$.

Recall that given a combinatorial data of an associative Belavin-Drinfeld structure $(S,C_1,C_2,A)$,
we have constructed a symplectic $2$-manifold $(\Sigma, \omega) $ and Lagrangians $L_1, L_2 \subset
\Sigma$. Recall also that, $\omega = \omega_g$ is the area form of a flat Riemannian metric $g$ on
$\Sigma$ and the Lagrangians $L_1, L_2$ are geodesic curves of length 1. 

\begin{defi}
    \label{pushoff} 
    Given $x$, $y \in \C$ we define the complex
    push-off $L_1^x$ of $L_1$ (resp. $L_2^y$ of $L_2$) to be the exact Lagrangian $L_1$ (resp.
    $L_2$) equipped with the complex rank 1 local system with monodromy $e^x$ (resp. $e^y$)
\end{defi}

Now, we let $\mathcal{X}$ to be the family of isomorphism classes of 
objects $\{ L_1^x \}$ for $x \in \C$, and similarly, we let $\mathcal{Y}$ to be the family of isomorphism classes of
objects $\{ L_2^y \}$ for $y \in \C$. For simplicity of notation, we sometimes write $x$ and $y$ for the
corresponding objects $L_1^x$ and $L_2^y$ of $\mathcal{F}(\Sigma)$. We remark that since by construction $L_1^x$ and $L_2^y$ are connected, gradable exact Lagrangians in $\Sigma$, up to shift there are unique grading structures on $L_1^x$ and $L_2^y$. We choose the shifts so that $CF(L_1^x, L_2^y)$ is supported in degree 0 for all $x,y$.

Note that we can apply the discussion from Section \ref{alg-form-AYBE-sec}
    to deduce that the family of objects in the Fukaya category over the formal disk,
  associated with an analytic family $(L_1^x)$ (resp.
    $(L_2^y)$) over $\C$, is
    quasi-isomorphic to the twisted object $(L_1, x \cdot \xi_{L_1})$ (resp. $(L_2, y \cdot \xi_{L_2})
    )$ as in Section \ref{twistedobjects}. (A sketch
    of a geometric proof of this also appears as \cite[Lemma 4.1]{auroux}.)

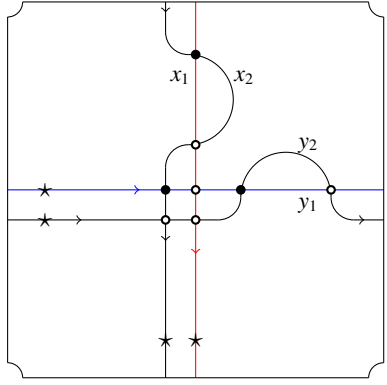
\begin{figure}[htb!]
\centering
\begin{tikzpicture} [scale=1]
	  \tikzset{->-/.style={decoration={ markings,
	    mark=at position #1 with {\arrow{>}}},postaction={decorate}}}

\draw[blue, ->-=.35]  (0,2.5) -- (5,2.5);

\draw [->-=.35] (0,2.1) -- (2.8,2.1);     
\draw [->-=.35] (4.6,2.1) -- (5,2.1);     
\draw (4.3,2.4) arc (0:180:0.6); 
\draw (3.1,2.4) arc (0:-90:0.3); 
\draw (4.3,2.4) arc (180:270:0.3); 

\draw [red, ->-=.67] (2.5,5) -- (2.5,0);
	
\draw [->-=.35] (2.1,2.8) -- (2.1,0);     
\draw [->-=.35] (2.1,5) -- (2.1,4.6);     
\draw (2.4,4.3) arc (90:-90:0.6); 
\draw (2.1,2.8) arc (180:90:0.3); 
\draw (2.1,4.6) arc (180:270:0.3);

\node at (2.5,0.5) {$\star$};
\node at (2.1,0.5) {$\star$};
\node at (0.5,2.1) {$\star$};
\node at (0.5,2.5) {$\star$};

\draw (0.2,0) -- (4.8,0);
\draw (0,0.2) -- (0,4.8);
\draw (0.2,5) -- (4.8,5);
\draw (5,0.2) -- (5,4.8);
\draw (0.2,0) arc (0:90:0.2); 
\draw (0,4.8) arc (-90:0:0.2);
\draw (4.8,5) arc (180:270:0.2);
\draw (4.8,0) arc (180:90:0.2);
\node at (3.15,4)   {\footnotesize $x_2$};
\node at (2.3,4)   {\footnotesize $x_1$};
\node at (4,3.1)   {\footnotesize $y_2$};
\node at (4, 2.3)   {\footnotesize $y_1$};

\draw[thick, fill=black] (2.5,4.3) circle(.05); 
\draw[thick, fill=white] (2.5,3.1) circle(.05); 
\draw[thick, fill=white] (4.3,2.5) circle(.05); 
\draw[thick, fill=black] (3.1,2.5) circle(.05); 
\draw[thick, fill=white] (2.5,2.5) circle(.05); 
\draw[thick, fill=black] (2.1,2.5) circle(.05); 
\draw[thick, fill=white] (2.5,2.1) circle(.05); 
\draw[thick, fill=white] (2.1,2.1) circle(.05);

\end{tikzpicture}
\caption{Hamiltonian perturbations of $L_1$ and $L_2$ (equipped with orientations and
    $\C^\times$-local systems)}
\label{fig3}
\end{figure}

Figure \ref{fig3} shows the simplest example on $\mathbb{T}_0$, where we have drawn four objects
$(x_1,x_2, y_1,y_2)$ in the punctured torus, which corresponds to Belavin-Drinfeld data with $S= \{
    1 \}, C_1 = (1), C_2=(1)$. Note that the underlying exact Lagrangians of $x_1$ and $x_2$ (resp.
$y_1$ and $y_2$) are Hamiltonian isotopic, however the monodromies of the complex local systems on
them are different.

Let us write $CF(x_1, x_2) = \C s_0 \oplus \C s_1$ and $CF(y_1,y_2) = \C t_0 \oplus \C t_1$. In
Figure \ref{fig3}, we denoted degree 0 generators by hollow and degree 1 generators by black dots
for these chain complexes.  In what follows, the signs come from the orientation of various moduli spaces, which we computed following the prescription in
    \cite[Sec. 7]{seidelgenus2}. 

We can compute the Floer differential to be:
\begin{align*} 
    \m_1( s_0) &= -s_1 + e^{x_2-x_1} s_1 \in CF^1(x_1,x_2)\\ 
    \m_1( t_0) &= -t_1 + e^{y_2-y_1} t_1 \in CF^1(y_1,y_2),
\end{align*}
where the terms correspond to the two visible lunes in each case. Hence, for $x_1 \neq x_2$ and $y_1 \neq y_2$, we have $HF(x_1,x_2) = HF(y_1,y_2)=0$. 

We also have $CF^{\neq 0} (x,y)=0$ for all $x,y$. Therefore, as explained in
the introduction, for distinct objects $x_1,x_2$ and $y_1,y_2$, the triple
Massey product: \[ \mathrm{MP} : CF^0(x_2,y_2) \otimes CF^1(y_1,x_2) \otimes CF^0(x_1,y_1)
\to CF^0(x_1,y_2) \] dualizes to a tensor 
\[ r^{x_1,x_2}_{y_1,y_2} : CF^0(x_2,y_2) \otimes CF^0(x_1,y_1) \to CF^0(x_1,y_2) \otimes
CF^0(x_2,y_1) \]  that satisfies
the AYBE.

Our next result (stated as Theorem B in the introduction) is that the obtained solution of the AYBE is precisely the trigonometric solution 
\eqref{trig-AYBE-eq} associated with $(S,C_1,C_2,A)$.

\begin{thm} Let $\Sigma$ be the square-tiled surface associated with an associative Belavin-Drinfeld structure $(S,C_1,C_2,A)$.
Then the tensor $r^{x_1,x_2}_{y_1,y_2}$ obtained from the triple products in the Fukaya category
$\mathcal{F}(\Sigma)$ only depends on $u=x_2-x_1$, $v=y_2-y_1$ and is a solution of the
AYBE over $\C$ given precisely by the formula \eqref{trig-AYBE-eq}.
\end{thm} 

{\it Proof.} The proof of this theorem follows from a direct computation of triple Massey products in the Fukaya
category $\mathcal{F}(\Sigma)$.

For clarity, we first do the computation for the simplest case, that is when
$S=\{p\}$ is a single point and $A$ is empty. Let us label the generators as follows:
    \[ CF^0(x_2,y_2) = \C \cdot p_{22} \ , \ CF^1(y_1,x_2)= \C \cdot q_{12} \ , \ CF^0(x_1,y_1) = \C \cdot
    p_{11} \ , \  CF^0(x_1,y_2) = \C \cdot p_{12} \]

Note that geometrically these generators correspond to the corners of the small square in the
    middle in Figure \ref{fig3}.  We are interested in computing the Massey product:
    \[ \mathrm{MP}(p_{22}, q_{12}, p_{11}) = \m_3 (p_{22},q_{12},p_{11}) - \m_2(h_2, p_{11}) -
    \m_2(p_{22}, h_1) \] 
where $h_1 \in CF^0(x_1,x_2)$ and $h_2 \in CF^0(y_2,y_1)$ satisfy $\m_1(h_1) =
    \m_2(q_{12},p_{11})$ and $\m_1(h_2) = \m_2(p_{22},q_{12})$. 

From Figure 3, it is straightforward to compute:
    \begin{align*}
        &\m_2(q_{12},p_{11}) = e^{x_2-x_1} \cdot s_1 \\
        &\m_2(p_{22},q_{12}) = t_1
    \end{align*} 
Therefore, we have
    \begin{align*}
        & h_1 = \left( \frac{e^{x_2-x_1}}{e^{x_2-x_1}-1} \right) \cdot s_0 \\
        & h_2 = \left( \frac{1}{e^{y_2-y_1}-1} \right) \cdot t_0 
    \end{align*} 
Again, from Figure 3, we can compute    
    \begin{align*}
        &\m_2(p_{22}, s_0) = p_{12}  \\ 
        &\m_2(t_0, p_{11}) = e^{y_2-y_1} \cdot p_{12} \\
        &\m_3(p_{22},q_{12},p_{11}) = p_{12} 
    \end{align*} 

Therefore, letting $u=x_2-x_1$, $v=y_2-y_1$, we conclude that
    \[ \mathrm{MP}(p_{22}, q_{12}, p_{11}) = \left( 1 + \frac{e^{u}}{1- e^u} + \frac{e^v}{1-e^v}
    \right)  \cdot p_{12} =  \left(\frac{1}{1- e^u} + \frac{1}{e^{-v}-1}
    \right) \cdot  p_{12}.  \]

Since the points $p_{ij}$, $i,j=1,2$ can be canonically identified with the intersection
    point $S= L_1 \cap L_2= \{ p \}$, we can view this tensor as
    \begin{align} \label{simplest} r^{x_1,x_2}_{y_1,y_2} &= \left( \frac{1}{\exp(u)-1} + \frac{1}{1-\exp(-v)} \right) e_{11}\ot e_{11} 
    \end{align}
(note that the dualization formula (\ref{dualiz}) brings in an extra overall sign).

In the case of general $S$ but with $A=\emptyset$, consider the covering map 
    \[ \pi : \Sigma_0 \to \mathbb{T}_0 \]
and the Lagrangians $X_i = \pi^{-1}(x_i)$ and $Y_i = \pi^{-1}(y_i)$ for $i=1,2$. 
Let us identify the points of intersection $X_i \cap Y_j$ by the set
$(\mathbf{e}_i)_{i \in S}$. Now, as before, we are interested in computing Massey products of the form 
\[ \mathrm{MP} : HF^0(X_2,Y_2) \otimes HF^1(Y_1,X_2) \otimes HF^0(X_1,Y_1) \to HF^0(X_1,Y_2) \]

Since Massey products are quasi-isomorphism invariants, we can compute each one with a convenient
Hamiltonian perturbation. Recall that we have the formula:
\[ \mathrm{MP}(\mathbf{e}_i, \mathbf{e}_j, \mathbf{e}_k) = \m_3(\mathbf{e}_i, \mathbf{e}_j, \mathbf{e}_k) -
\m_2( h_2, \mathbf{e}_k) - \m_2(\mathbf{e}_i, h_1) \]
where $\m_1(h_2) = \m_2 (\mathbf{e}_i, \mathbf{e}_j)$ and $\m_2(h_1) =
\m_2(\mathbf{e}_j, \mathbf{e}_k)$. 

We first observe that $\m_2(\mathbf{e}_i, \mathbf{e}_j)=0$ if $i \neq j$ since there are no triangles
that can contribute by construction, and
$\m_3(\mathbf{e}_i, \mathbf{e}_j, \mathbf{e}_k) =0$ unless $i=j$ or $j=k$, since $A$ is empty.

Therefore, the only possibly non-trivial triple Massey products are of the form:
\begin{align*} &\mathrm{MP}(C_1^k(\mathbf{e}_i) , \mathbf{e}_i , \mathbf{e}_i)  \text{\
    \ \ for } k=0,1,\ldots (n-1) \\
    &\mathrm{MP}(\mathbf{e}_i, \mathbf{e}_i, C_2^m(\mathbf{e}_i)) \text{\ \ \  for } m=0,1\ldots
    (n-1). 
\end{align*} 

For ease of computation, we arrange that the holonomy contributions of the
$\C^\times$-local systems on $L_i$ are divided equally to $n$ parts, each contributing $e^{\frac{u}{n}}$ for $L_1$ and $e^{\frac{v}{n}}$ for $L_2$, interlaced between the
$n$ intersection points $L_i \cap L_j$. In other words, each time a holomorphic disk has
boundary covering one of these regions, there is an associated weight $e^{\pm \frac{u}{n}}$ or
$e^{\pm \frac{v}{n}}$, where the sign of the exponent is determined, as before, according to whether the boundary orientation of the holomorphic disk matches that of $L_i$ or not. 

The computation of $\mathrm{MP}(\mathbf{e}_i, \mathbf{e}_i, \mathbf{e}_i)$ is done in a completely
analogous way to the above computation given for $n=1$, hence we have:
\[ \mathrm{MP}(\mathbf{e}_i, \mathbf{e}_i, \mathbf{e}_i) =  \left(\frac{1}{1- e^u} +
\frac{1}{e^{-v}-1}  \right) \mathbf{e}_i \]

\begin{figure}[htb!]
\centering
\begin{tikzpicture} [scale=1]
	  \tikzset{->-/.style={decoration={ markings,
	    mark=at position #1 with {\arrow{>}}},postaction={decorate}}}
	
\filldraw[fill=blue!40!white, fill opacity=0.4, draw opacity=0] (1.2,1.2) rectangle (5,2.2);

\filldraw[fill=blue!40!white, fill opacity=0.4, draw opacity=0] (0,3.2) rectangle (2.2,4.2);

\draw [->-=.35] (1.2,2.2) -- (5,2.2);     
\draw [->-=.35] (1.2,1.2) -- (5,1.2);	
\draw [->-=.35] (0,4.2) -- (2.2,4.2);     
\draw [->-=.35] (0,3.2) -- (2.2,3.2);	
    
\draw [->-=.67] (2.2,5) -- (2.2,0);
\draw [->-=.67] (1.2,5) -- (1.2,0);

\node at (1,5)   {\footnotesize $x_2$};
\node at (2,5)   {\footnotesize $x_1$};
\node at (4, 2.4)   {\footnotesize $y_1$};
\node at (4, 1)   {\footnotesize $y_2$};
\node at (0.5, 4.4)   {\footnotesize $y_1$};
\node at (0.5, 3)   {\footnotesize $y_2$};

\node at (2.6,3.2) {\tiny $C_2^m(\mathbf{e}_i)$}; 
\node at (2.6,4.2) {\tiny $C_2^m(\mathbf{e}_i)$}; 
\node at (1,2.2) {\tiny $\mathbf{e}_i$}; 
\node at (1,1.2) {\tiny $\mathbf{e}_i$};

\filldraw[fill=blue!40!white, fill opacity=0.4, draw opacity=0] (8.2,1.2) rectangle (9.2,4.2);

\draw [->-=.35] (8.2,2.2) -- (12.2,2.2);     
\draw [->-=.35] (8.2,1.2) -- (12.2,1.2);	
\draw [->-=.35] (7.2,4.2) -- (9.2,4.2);     
\draw [->-=.35] (7.2,3.2) -- (9.2,3.2);	
\
    
\draw [->-=.67] (9.2,5) -- (9.2,0);
\draw [->-=.67] (8.2,5) -- (8.2,0);

\node at (8,5)   {\footnotesize $x_2$};
\node at (9,5)   {\footnotesize $x_1$};
\node at (11, 2.4)   {\footnotesize $y_1$};
\node at (11, 1)   {\footnotesize $y_2$};
\node at (7.5, 4.4)   {\footnotesize $y_1$};
\node at (7.5, 3)   {\footnotesize $y_2$};

\node at (8.35,4.33) {\tiny $\mathbf{e}_i$}; 
\node at (9.35,4.33) {\tiny $\mathbf{e}_i$}; 
    \node at (9.55,1) {\tiny $C_1^k(\mathbf{e}_i)$}; 
    \node at (8.55,1) {\tiny $C_1^k(\mathbf{e}_i)$}; 

\end{tikzpicture}
\caption{Horizontally (left) and vertically (right) extending rectangles}
\label{fig4}
\end{figure}
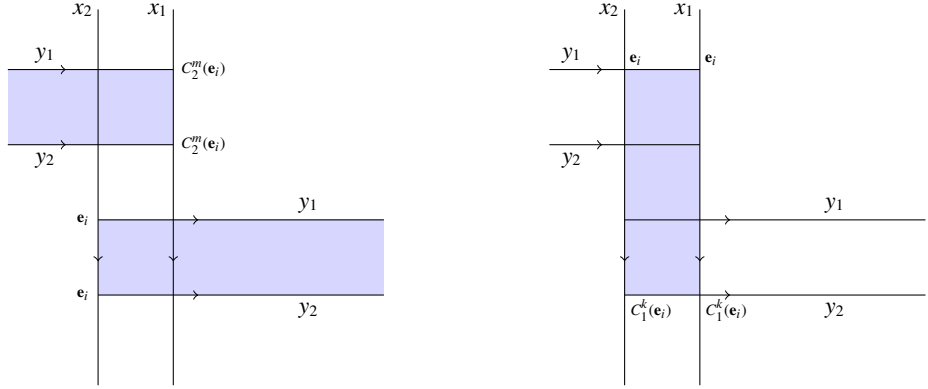

Next, we observe that there are two families of rectangles with boundary on
$(X_1,Y_2,X_2,Y_1)$ as illustrated in Figure \ref{fig4}. 

These contribute to $\m_3$ products of the form
\[ \m_3( C_1^k(\mathbf{e}_i) , \mathbf{e}_i, \mathbf{e}_i) = e^{\frac{ku}{n}}
C_1^k(\mathbf{e}_i)  \]
for $i =1,\ldots, n$ and $k=1,\ldots, n-1$ and 
\[ \m_3(\mathbf{e}_i, \mathbf{e}_i, C_2^m(\mathbf{e}_i))
=e^{\frac{mv}{n}} C_2^m(\mathbf{e}_i)  \]
for $i=1,\ldots, n$ and $m=1,\ldots, n-1$. 

Furthermore, as before let $CF(X_1,X_2)= \C s_0 \oplus \C s_1$ and $CF^0(Y_1,Y_2) = \C t_0 \oplus \C t_1 $. We compute the products:
\begin{align*}
    \m_2 : CF^1(Y_1,X_2) \otimes CF^0(X_1,Y_1) &\to CF^0(X_1,X_2) \\
    \m_2( \mathbf{e}_i, \mathbf{e}_i ) &= e^{\frac{ku}{n}} s_1, \\
    \m_2 : CF^0(X_2,Y_2) \otimes CF^1(Y_1,X_2) &\to CF^0(Y_1,Y_2) \\
    \m_2( \mathbf{e}_i, \mathbf{e}_i) &= t_1, 
\end{align*}
and 
\begin{align*}
    &\m_2(\mathbf{e}_{i}, s_0) = \mathbf{e}_1, \\ 
    &\m_2(t_0, C_2^{m} (\mathbf{e}_{i})) = e^{\frac{mv}{n}} \mathbf{e}_1
\end{align*}
Thus, we conclude that 
\[ \mathrm{MP}( \mathbf{e}_i , C_1^k(\mathbf{e}_i), C_1^k(\mathbf{e}_i)) = e^{\frac{ku}{n}}
 \left(1 + \frac{e^u}{1-e^{u}} \right) \mathbf{e}_i  \]
for $i =1,\ldots, n$ and $k=1,\ldots, n-1$ and 
\[ \mathrm{MP}(\mathbf{e}_i, \mathbf{e}_i, C_2^m(\mathbf{e}_i))
=e^{\frac{mv}{n}} \left( 1+\frac{e^v}{1-e^v} \right) C_2^m(\mathbf{e}_i )    \]
for $i=1,\ldots, n$ and $m=1,\ldots, n-1$. 

Dualising to the tensor $r^{x_1,x_2}_{y_1,y_2}$, we get the terms:
\begin{comment}
\begin{equation}   
\begin{aligned}
    &q^{\mathrm{Re}(u)\mathrm{Re}(v)} \biggl( \sum_{i} (1 + \frac{1}{1-\exp(-u)} + \frac{1}{1-\exp(-v)} )
    e_{ii} \otimes e_{ii} \biggr. \\ 
    &+ \frac{1}{1-\exp(-u)} \sum_{0<k<n, i}\biggl.  \exp(-\frac{ku}{n})e_{ii}\ot e_{C_1^k(i),C_1^k(i)} +
      \frac{1}{1-\exp(-v)}\sum_{0<m<n, i}\exp(-\frac{mv}{n})e_{C_2^m(i),i}\ot e_{i, C_2^m(i)} \biggr)
    \end{aligned}
    \end{equation}
Noting that $1 + \frac{1}{1-\exp(-u)} = \frac{1}{\exp(u)-1}$, and letting $i \to C_1^{n-k}(i)$, $k \to
n-k$ and  
$i \to C_2^{m-k}(i)$, $m \to n-m$ in the last two sums respectively, we can rewrite this as:
\end{comment} 
\begin{equation}   
\begin{aligned}
    &\sum_{i} (\frac{1}{\exp(u)-1} + \frac{1}{1-\exp(-v)} )
    e_{ii} \otimes e_{ii} \biggr. \\ 
    &+ \frac{1}{\exp(u)-1} \sum_{0<k<n, i}\biggl.  \exp(\frac{ku}{n}) e_{C_1^k(i),C_1^k(i)} \ot
    e_{ii} +
      \frac{1}{\exp(v)-1}\sum_{0<m<n, i}\exp(\frac{mv}{n}) e_{i, C_2^m(i)} \ot e_{C_2^m(i),i} \end{aligned}
    \end{equation}
and we see that this agrees with the stated result in the case $A$ is empty.  

\begin{figure}[htb!]
\centering
\begin{tikzpicture} [scale=1]
	  \tikzset{->-/.style={decoration={ markings,
	    mark=at position #1 with {\arrow{>}}},postaction={decorate}}}
	
\filldraw[fill=blue!40!white, fill opacity=0.4, draw opacity=0] (0.8,4.2) rectangle (4.2,0.8);

\draw [->-=.27] (0,1.8) -- (5,1.8);     
\draw [->-=.27] (0,0.8) -- (5,0.8);	
\draw [->-=.27] (0,4.2) -- (5,4.2);     
\draw [->-=.27] (0,3.2) -- (5,3.2);

\draw [->-=.75] (4.2,5) -- (4.2,0);
\draw [->-=.75] (3.2,5) -- (3.2,0);
\draw [->-=.75] (1.8,5) -- (1.8,0);
\draw [->-=.75] (0.8,5) -- (0.8,0);

\node at (0.6,5)   {\footnotesize $x_2$};
\node at (2,5)   {\footnotesize $x_1$};

\node at (3,0)   {\footnotesize $x_2$};
\node at (4.4,0)   {\footnotesize $x_1$};

\node at (4.8, 2)   {\footnotesize $y_1$};
\node at (4.8, 0.6)   {\footnotesize $y_2$};
\node at (0.2, 4.4)   {\footnotesize $y_1$};
\node at (0.2, 3)   {\footnotesize $y_2$};

\node at (4.7,1) {\tiny $C_1^k C_2^m(a)$}; 
    \node at (4.55,4) {\tiny $C_2^m(a)$}; 
    \node at (0.65, 4) {\tiny $a$}; 
    \node at (0.45,0.95) {\tiny $C_1^k(a)$}; 

\filldraw[fill=blue!40!white, fill opacity=0.4, draw opacity=0] (8.8,3.2) rectangle (10.2,1.8);

\draw [->-=.27] (7,1.8) -- (12,1.8);     
\draw [->-=.27] (7,0.8) -- (12,0.8);	
\draw [->-=.27] (7,4.2) -- (12,4.2);     
\draw [->-=.27] (7,3.2) -- (12,3.2);

\draw [->-=.75] (11.2,5) -- (11.2,0);
\draw [->-=.75] (10.2,5) -- (10.2,0);
\draw [->-=.75] (8.8,5) -- (8.8,0);
\draw [->-=.75] (7.8,5) -- (7.8,0);

\node at (7.6,5)   {\footnotesize $x_2$};
\node at (9,5)   {\footnotesize $x_1$};

\node at (10,0)   {\footnotesize $x_2$};
\node at (11.4,0)   {\footnotesize $x_1$};

\node at (11.8, 2)   {\footnotesize $y_1$};
\node at (11.8, 0.6)   {\footnotesize $y_2$};
\node at (7.2, 4.4)   {\footnotesize $y_1$};
\node at (7.2, 3)   {\footnotesize $y_2$};

\node at (10.7,2) {\tiny $C_1^k C_2^m(a)$}; 
    \node at (10.55,3) {\tiny $C_2^m(a)$}; 
    \node at (8.65, 3) {\tiny $a$}; 
\node at (8.45,1.95) {\tiny $C_1^k(a)$};

\end{tikzpicture}
\caption{The two rectangles for each $a \in A(m,k)$ }
\label{fig5}
\end{figure}
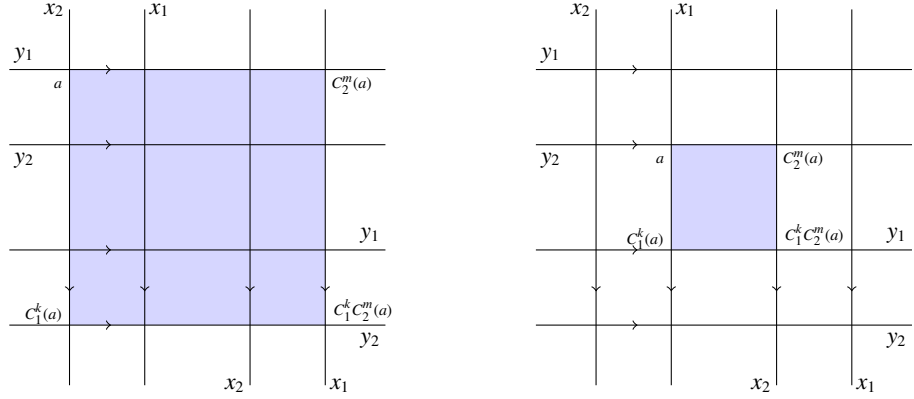 

Finally, we will compute the contribution of rectangles when $A$ is non-empty. When $A$ is
non-empty, the rectangular regions with corners $(a, C_1(a), C_1C_2(a), C_2(a))$ are filled for each
$a \in A$. We get two new contributions to the $\m_3$ product from such regions (See Figure
\ref{fig5}). Furthermore, it
may happen that union of those regions also give new rectangles. A combinatorial way to encode this
is to let $A(k,m) \subset A$ to be the set of all $a \in A$ such that $C_1^i C_2^j(a)
\in A$ for all $0 \leq i < k$, $0 \leq j <m$, then for each $a \in A(k,m)$ we have the following
contributions because of the filled rectangular region with corners $(a, C_1^k(a), C_1^k
C_2^m (a), C_2^m(a))$:
\[ \m_3( C_1^{k}(a) , a, C_2^m(a)) = 
e^{\frac{ku  + mv}{n}} C_1^{k} C_2^m(a) \]
corresponding to the rectangle drawn on the left of Figure \ref{fig5}, and  
\[ \m_3( C_2^m(a) , C_1^k C_2^m(a), C_1^{k}(a)) = -e^{-\frac{ku+ mv}{n}} a  \]
corresponding to the rectangle drawn on the right of Figure \ref{fig5}. 

The signs that appear in the two formulae are affected by the orientations of
the Lagrangians and we note that unlike the appearance, there is no typographical error in what we
wrote. The sign in the exponentials are determined according to whether the orientation of the
Lagrangians agree with the counter-clockwise boundary orientation of the rectangle, and the overall
sign is determined according to the orientation of the moduli space which we computed as before
using \cite[Sec. 7]{seidelgenus2}. 

Recall also that the dualization formula (\ref{dualiz}) brings in an extra overall sign. Thus, we conclude that in the case of arbitrary $A$ we have in addition the contribution of the
following terms to $r(u,v)$, indexed by elements of the subsets $A(k,m)$, $k,m>0$: 
\begin{align}
    \sum_{\substack{0<k,0<m; \\ a\in A(k,m)}} \Big\{ \exp(-\frac{ku+mv}{n})e_{C_2^m(a),a}\ot
e_{C_1^k(a),C_1^k C_2^m(a)} -\exp(\frac{ku+mv}{n})e_{C_1^k(a),C_1^kC_2^m(a)}\ot e_{C_2^m(a),a} \Big\}        
\end{align}
\ed

\begin{rem} We would like to mention an alternative to the above computation. It may appear more
    natural to take complex push-offs of $L_1$ and $L_2$ as follows.  First, on $\mathbb{T}_0$, let
    $l_1^x$ (resp. $l_2^y$)  be the geodesic push-off of $l_1$ (resp. $l_2$) such that the oriented area of the cylinder
    bounded by $l_1$ and $l_1^x$ (resp. $l_2$ and $l_2^y$) is $\Re(x)$ (resp. $\Re(y)$). We then set
    $L_1^x =
    \pi^{-1}(l_1^x)$ (resp. $L_2^y = \pi^{-1}(l_2^y)$) equipped with a $\mathrm{U}(1)$-local system with
monodromy $e^{i \Im(x)}$ (resp. $e^{i \Im(y)}$). The simplest case is shown in
    Figure \ref{fig6}. 

\begin{figure}[htb!]
\centering
\begin{tikzpicture} [scale=1]

    \tikzset{->-/.style={decoration={ markings,
	    mark=at position #1 with {\arrow{>}}},postaction={decorate}}}
\draw (0.2,0) -- (4.8,0);
\draw[blue, ->-=.35]  (0,2.5) -- (5,2.5);
\draw (0,0.2) -- (0,4.8);
\draw [red, ->-=.67] (2.5,5) -- (2.5,0);
	
\draw [->-=.35] (0,2.2) -- (5,2.2);     
\draw [->-=.35] (0,1.2) -- (5,1.2);	
\draw [->-=.67] (2.2,5) -- (2.2,0);
\draw [->-=.67] (1.2,5) -- (1.2,0);

\node at (2.205,0.5) {$\star$};
\node at (1.205,0.5) {$\star$};
\node at (0.5,2.195) {$\star$};
\node at (0.5,1.195) {$\star$};

\draw (0.2,5) -- (4.8,5);
\draw (5,0.2) -- (5,4.8);
\draw (0.2,0) arc (0:90:0.2); 
\draw (0,4.8) arc (-90:0:0.2);
\draw (4.8,5) arc (180:270:0.2);
\draw (4.8,0) arc (180:90:0.2);
\node at (2.8,4)   {\footnotesize $L_1$};
\node at (1,4)   {\footnotesize $x_2$};
\node at (2,4)   {\footnotesize $x_1$};
\node at (4,2.75)   {\footnotesize $L_2$};
\node at (4, 2)   {\footnotesize $y_1$};
\node at (4, 1)   {\footnotesize $y_2$};

\end{tikzpicture}
\caption{Non-exact push-offs of $L_1$ and $L_2$ (equipped with orientations and $\mathrm{U}(1)$-local systems)}
\label{fig6}
\end{figure}
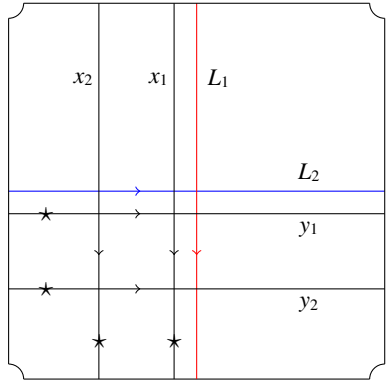 

However, there is a significant drawback in this approach. Namely, the push-offs $L_1^x$ and
    $L_1^y$ are no longer exact Lagrangians when $\Re(x), \Re(y) \neq 0$. Hence, we cannot guarantee a
    priori that the count of holomorphic disks is finite (or convergent). Therefore, in
    this set-up one has to work over the Novikov field 
    \[  \Lambda = \{ \sum_{i \in \Z} a_i q^{t_i} |a_i \in \C, a_i = 0 \text{ for } i \ll 0, t_i
    \in \R, t_i \to \infty \}  \]
    and the formula (\ref{matrixcoef}) should be modified as
\begin{equation} \m_k(\rho_k,\ldots,\rho_1) = \sum_{\substack{ [u]: \text{ind}([u])=2-k}  } 
\# \mathcal{M}(p_k,p_{k-1},\ldots,p_1,p_0; [u]) \cdot \text{hol}_{\partial u} \cdot q^{\int_u \omega}. \end{equation}   
With this in place, one can compute the corresponding Massey product simply by counting rectangles. In the simplest case, that is when $S = \{p \}$ is a single point and $A$ is
    empty, computing the tensor $r^{x_1, x_2}_{y_1,y_2}$ boils down to counting rectangles with
    corners $(p_{12} , p_{22} , q_{12} , p_{11})$  in the counter-clockwise order
      weighted by their
      areas. Interestingly, there are indeed infinitely many homotopy classes of rectangles that contribute to this count. The smallest rectangle with corners in $(p_{12}, p_{22}, q_{12}, p_{11})$ has area
      $\mathrm{Re}(u)\mathrm{Re}(v)$. Other than this, there are two families of rectangles - namely, those that are
      horizontally or
      vertically extending. Writing $x_i= a_i + i \alpha_i$ and $y_i = b_i + i \beta_i$, the horizontally extending ones
    are weighted by $e^{il(\alpha_2-\alpha_1)} q^{l(a_1-a_2) + (a_1-a_2)(b_1-b_2)}$ for
    $l=1,2,\ldots$, and the vertically extending ones are weighted by $e^{il(\beta_2-\beta_1)} q^{l(b_1-b_2)
    + (a_1-a_2)(b_1-b_2)}$ for $l=1,2,\ldots.$

The overall contribution of all these rectangles can be computed as: 
    \[ r^{x_1,x_2}_{y_1,y_2} (p_{11} \otimes p_{22}) =
     -q^{(a_1-a_2)(b_1-b_2)} \left( 1 + \sum_{l=1}^\infty
    e^{li(\alpha_2-\alpha_1)} q^{l(a_1-a_2)} + \sum_{l=1}^\infty e^{li(\beta_2-\beta_1)}
    q^{l(b_1-b_2)} \right) (p_{12} \otimes p_{21}) \]

Since the points $p_{ij}$, $i,j=1,2$ can be canonically identified with the intersection
    point $S= L_1 \cap L_2= \{ p \}$, we can view this tensor as
    \begin{align*} r^{x_1,x_2}_{y_1,y_2} &= -q^{(a_1-a_2)(b_1-b_2)} \left( 1 + \sum_{l=1}^\infty
        e^{li(\alpha_2-\alpha_1)}
    q^{l(a_1-a_2)} + \sum_{l=1}^\infty e^{li(\beta_2-\beta_1)} q^{l(b_1-b_2)} \right) (e_{11}
        \otimes e_{11}) 
\end{align*}
Now, we observe that the series expansion has positive radius of convergence equal to 1, hence in particular specializing the Novikov parameter $q =e^{-1}$ leads to the formula:
\[  r^{x_1,x_2}_{y_1,y_2} = -e^{-\mathrm{Re}(u) \mathrm{Re}(v)}  \left( 1 + \frac{e^{u}}{1-e^{u}} + \frac{e^v}{1-e^v} \right) e_{11}\ot e_{11}  = e^{-\mathrm{Re}(u) \mathrm{Re}(v)}  \left( \frac{1}{1-e^{u}} + \frac{1}{e^{-v}-1} \right) e_{11}\ot e_{11}    \]
This is remarkably in agreement with what we have computed before in formula (\ref{simplest}) up to the overall constant $e^{-\mathrm{Re}(u) \mathrm{Re}(v)}$ which can be absorbed into the choice of basis. Similar computation can be done in all cases. This gives a hint that in an appropriately defined Fukaya category, the two different ways of pushing-off $L_1$ and $L_2$ should lead to quasi-isomorphic objects. (Compare with the discussion in \cite[Section 4.1]{auroux}.) 
\end{rem}

\begin{rem} Since the $A_\infty$-relations hold in the Fukaya category by its general construction, Theorem $B$ gives a new proof of the fact that $r(u,v)$ given by \eqref{trig-AYBE-eq} satisfies
the AYBE, which is proved in \cite{pol09} by a rather tedious calculation.
 On the other hand, in \cite{pol09} it was also proven that for $r(u,v)$ given by \eqref{trig-AYBE-eq},
 %is a   strongly nondegenerate, trigonometric solution to AYBE (as given in Theorem B by an explicit formula), then 
    \[ R(u,v) = \left( \frac{(e^{\frac{u}{2}}- e^{\frac{-u}{2}}) \cdot (e^{\frac{v}{2}}-
    e^{\frac{-v}{2}})}{ e^{\frac{u}{2}}- e^{\frac{-u}{2}} + e^{\frac{v}{2}}-
    e^{\frac{-v}{2}}} \right) r(u,v) \] 
    satisfies the quantum Yang-Baxter equation (for fixed $u$):
    \[ R^{12}(v) R^{13}(v+v') R^{23}(v') = R^{23}(v') R^{13}(v+v') R^{12}(v) \]
with the unitarity condition
    \[ R(u,v) R^{21}(u, -v) = 1 \otimes 1. \]
We do not know a conceptual explanation for this. It would be interesting to study this in the
    setting of Fukaya categories. 
\end{rem}

We will need the following result in the proof of Theorem C.

\begin{prop}\label{g1-Fuk-gen-prop}
Let $(S,C_1,C_2,A)$ be an associative Belavin-Drinfeld structure, such that $C_1$ and $C_2$ commute, and let
$(\Sigma,L_1,L_2)$ be the corresponding square-tiled surface with a pair of Lagrangians (where $\Sigma$ is a punctured
    torus). Then $(L_1,L_2)$ split generates the Fukaya category
    $\mathcal{F}(\Sigma)$ of exact, compact (graded) Lagrangians in $\Sigma$.
\end{prop}

\begin{figure}[htb!]
\centering
\begin{tikzpicture} [scale=1.5]

        \tikzset{->-/.style={decoration={ markings,
	    mark=at position #1 with {\arrow{>}}},postaction={decorate}}}
        \tikzset{->>-/.style={decoration={ markings,
	    mark=at position #1 with {\arrow{>>}}},postaction={decorate}}}
        \tikzset{->>>-/.style={decoration={ markings,
	    mark=at position #1 with {\arrow{>>>}}},postaction={decorate}}}
       \tikzset{->>>>-/.style={decoration={ markings,
	    mark=at position #1 with {\arrow{>>>>}}},postaction={decorate}}}
        \tikzset{->>>>>-/.style={decoration={ markings,
	    mark=at position #1 with {\arrow{>>>>>}}},postaction={decorate}}}

       \tikzset{-|>-/.style={decoration={ markings,
	    mark=at position #1 with {\arrow{latex}}},postaction={decorate}}}

	\draw[thick, ->-=.55] (-15,0) -- (-14,0);
    \draw[thick, ->>-=.6] (-14,0) -- (-13,0);
    \draw[thick, ->>>-=.65] (-13,0) -- (-12,0);
    \draw[thick, ->>>>-=.7] (-12,0) -- (-11,0);
    \draw[thick, ->>>>>-=.7] (-11,0) -- (-10,0);
        
    \draw[thick, ->>>>-=.65] (-15,1) -- (-14,1);
    \draw[thick, ->>>>>-=.7] (-14,1) -- (-13,1);
    \draw[thick, ->-=.6] (-13,1) -- (-12,1);
    \draw[thick, ->>-=.55] (-12,1) -- (-11,1);
    \draw[thick, ->>>-=.55] (-11,1) -- (-10,1);

	\draw[thick, -|>-=.6] (-15,0) -- (-15,1);
	\draw[thick, -|>-=.6] (-10,0) -- (-10,1);

	\draw[thick] (-14,0) -- (-14,1);
	\draw[thick] (-13,0) -- (-13,1); 
	\draw[thick] (-12,0) -- (-12,1); 
	\draw[thick] (-11,0) -- (-11,1); 
         
	\draw[red, ->-=.7] (-14.5,1) -- (-14.5,0); 
    \draw[red, ->-=.7] (-13.5,1) -- (-13.5,0); 
	\draw[red, ->-=.7] (-12.5,1) -- (-12.5,0); 
    \draw[red, ->-=.7] (-11.5,1) -- (-11.5,0); 
    \draw[red, ->-=.7] (-10.5,1) -- (-10.5,0); 

    \draw[green!50!black, ->-=.7] (-14.5,1) -- (-15,0.75); 
    \draw[green!50!black, ->-=.7] (-13.5,1) -- (-15,0.25); 
	\draw[green!50!black, ->-=.7] (-12.5,1) -- (-14.5,0); 
    \draw[green!50!black, ->-=.7] (-11.5,1) -- (-13.5,0); 
    \draw[green!50!black, ->-=.7] (-10.5,1) -- (-12.5,0); 
    \draw[green!50!black, ->-=.7] (-10,0.75) -- (-11.5,0); 
    \draw[green!50!black, ->-=.7] (-10,0.25) -- (-10.5,0);

    \draw[blue, ->-=.15, ->-=.35, ->-=.55, ->-=.75, ->-=.95] (-15,0.5) -- (-10,0.5); 

    \node[blue] at (-15.2,0.5)   {\footnotesize $L_1$};
    \node[red] at (-14.5,1.2)   {\footnotesize $L_2$};

\node[green!50!black] at (-14.5,-0.2)   {\footnotesize $M_1$};
\node[green!50!black] at (-13.5,-0.2)   {\footnotesize $M_2$};
\node[green!50!black] at (-12.5,-0.2)   {\footnotesize $M_3$};
\node[green!50!black] at (-11.5,-0.2)   {\footnotesize $M_4$};
\node[green!50!black] at (-10.5,-0.2)   {\footnotesize $M_5$};

        \draw[thick, fill=white] (-15,0) circle(.1); 
        \draw[thick, fill=white] (-14,0) circle(.1); 
        \draw[thick, fill=white] (-13,0) circle(.1); 
        \draw[thick, fill=white] (-12,0) circle(.1); 
        \draw[thick, fill=white] (-11,0) circle(.1); 
        \draw[thick, fill=white] (-10,0) circle(.1); 
        \draw[thick, fill=white] (-15,1) circle(.1); 
        \draw[thick, fill=white] (-14,1) circle(.1); 
        \draw[thick, fill=white] (-13,1) circle(.1); 
        \draw[thick, fill=white] (-12,1) circle(.1); 
        \draw[thick, fill=white] (-11,1) circle(.1); 
        \draw[thick, fill=white] (-10,1) circle(.1); 
        
        \node at (-15,0) {\tiny $A$}; 
        \node at (-14,0) {\tiny $B$}; 
        \node at (-13,0) {\tiny $C$};
        \node at (-12,0) {\tiny $D$}; 
        \node at (-11,0) {\tiny$E$}; 
        \node at (-10,0) {\tiny $A$}; 
        \node at (-15,1) {\tiny $D$}; 
        \node at (-14,1) {\tiny $E$}; 
        \node at (-13,1) {\tiny $A$}; 
        \node at (-12,1) {\tiny $B$}; 
        \node at (-11,1) {\tiny $C$}; 
        \node at (-10,1) {\tiny $D$}; 
    
\end{tikzpicture}
    \caption{Generators of the exact Fukaya category, $(n,k) = (5,2)$}
	\label{fig7}
\end{figure}
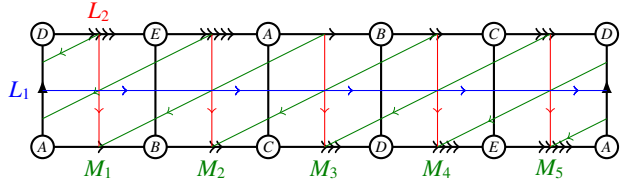

\Pf . We first prove that $(L_1,L_2)$ split generate when $A =\emptyset$ and $\Sigma =
\Sigma_0$. Without loss of generality, suppose that $S= \{1,\ldots, n\}$, $C_1(i)=i+1$ and that $C_2 = C_1^k$
for some $k$ which is prime to $n$. We can draw the corresponding square-tiled surface as in Figure
\ref{fig7} (where the case of $n=5$ and $k=2$ is drawn). Let $M_1,M_2,\ldots, M_n$ be $n$ disjoint
Lagrangians corresponding to curves of slope $1/k$, drawn in green in Figure \ref{fig7}. 

Note that these Lagrangians have a natural grading structure (since our line field is given by the horizontal foliation).
It was proven in \cite[Lem.\ 3.1.1]{LP} that the collection $L_1, M_1, M_2,\ldots, M_n$ split generates the exact Fukaya category. (This essentially follows from the fact that Dehn twists around $L_1, M_1, M_2,\ldots, M_n$ generate the pure mapping class group of the $n$-punctured torus.)

Note that we have the following intersections in homology:
\begin{align*}
    [M_i] \cdot [L_1] &= 1 \\
    [M_i] \cdot [L_2] &= k \\
    [L_2] \cdot [L_1] &= n
\end{align*}
In fact, by considering dual curves to $M_i$, it is easy to see that
\[ [L_2] = k [L_1] + [M_1] + [M_2] + \cdots + [M_n]  \in H_1(\Sigma_0) \] 

We claim that there is an exact triangle of the form:
\begin{equation} \label{eq:exactness}
    \xymatrix{
        M_1 \oplus M_2 \oplus \ldots \oplus M_n \ar[r] & L_1^{\oplus k} \ar[d] \\
                                & L_2 \ar[ul]^{[1]}
                                    } 
\end{equation}
where the maps $M_i \to L_1^{\oplus k}$ are given by $(c_i,c_i,\ldots, c_i)$, with $c_i \in
CF^{1}(M_i,L_1)$, for each $i$, being the generator corresponding to the unique intersection point. It is then clear
the $L_1$ and $L_2$ split generate $\mathcal{F}(\Sigma)$. 

The exact triangle is an example of a surgery exact triangle proven in this case by Abouzaid in
\cite[Lemma 5.4]{abouz}. It is technically easier to show that the following equivalent statement
holds:
\[ Cone( \ldots Cone((Cone( M_1 \oplus M_2 \oplus \ldots \oplus M_n) \to L_1) \to L_1)  \ldots \to
L_1 ) \simeq L_2. \]
Indeed, we first do a surgery at each intersection point of $L_1$ and each $M_i$ and
then we perform a new surgery at the $n$ intersection points of the obtained Lagrangian with a new
copy of $L_1$. We do this $k$ times (including the first surgery between $L_1$ and $M_i$'s) until we arrive at an exact Lagrangian Hamiltonian isotopic to $L_2$. Note that in each isotopy class of homotopically essential (i.e. not null-homotopic) simple closed curves, there is a unique exact Lagrangian up to Hamiltonian isotopy, so it suffices to
check that the end result of all the surgeries, which is an exact Lagrangian, is smoothly isotopic
to $L_2$. 

The corresponding picture is drawn in Figures \ref{fig8} and \ref{fig9} below for $n=5, k=2$ case, from which it is clear how the general
case works.

Note that by Seidel's exact triangle \cite{seidelLES}, the first iteration can be identified as 
\[Cone((M_1 \oplus M_2 \oplus \ldots \oplus M_n) \to L_1) \simeq \tau_{M_1} \circ \tau_{M_2} \circ
\ldots \circ \tau_{M_n}
    (L_1). \]

\begin{figure}[htb!]
\centering
\begin{tikzpicture} [scale=1.5]

        \tikzset{->-/.style={decoration={ markings,
	    mark=at position #1 with {\arrow{>}}},postaction={decorate}}}
        \tikzset{->>-/.style={decoration={ markings,
	    mark=at position #1 with {\arrow{>>}}},postaction={decorate}}}
        \tikzset{->>>-/.style={decoration={ markings,
	    mark=at position #1 with {\arrow{>>>}}},postaction={decorate}}}
       \tikzset{->>>>-/.style={decoration={ markings,
	    mark=at position #1 with {\arrow{>>>>}}},postaction={decorate}}}
        \tikzset{->>>>>-/.style={decoration={ markings,
	    mark=at position #1 with {\arrow{>>>>>}}},postaction={decorate}}}

       \tikzset{-|>-/.style={decoration={ markings,
	    mark=at position #1 with {\arrow{latex}}},postaction={decorate}}}

	\draw[thick, ->-=.55] (-15,0) -- (-14,0);
    \draw[thick, ->>-=.6] (-14,0) -- (-13,0);
    \draw[thick, ->>>-=.65] (-13,0) -- (-12,0);
    \draw[thick, ->>>>-=.7] (-12,0) -- (-11,0);
    \draw[thick, ->>>>>-=.7] (-11,0) -- (-10,0);
        
    \draw[thick, ->>>>-=.65] (-15,1) -- (-14,1);
    \draw[thick, ->>>>>-=.7] (-14,1) -- (-13,1);
    \draw[thick, ->-=.6] (-13,1) -- (-12,1);
    \draw[thick, ->>-=.55] (-12,1) -- (-11,1);
    \draw[thick, ->>>-=.55] (-11,1) -- (-10,1);

	\draw[thick, -|>-=.6] (-15,0) -- (-15,1);
	\draw[thick, -|>-=.6] (-10,0) -- (-10,1);

	\draw[thick] (-14,0) -- (-14,1);
	\draw[thick] (-13,0) -- (-13,1); 
	\draw[thick] (-12,0) -- (-12,1); 
	\draw[thick] (-11,0) -- (-11,1);

\draw[green!50!black, ->-=.7] (-14.5,1) -- (-15,0.75); 
\draw[green!50!black, ->-=.7] (-14.75,0.375) -- (-15,0.25); 
\draw[green!50!black, ->-=.7] (-13.5,1) -- (-14.25,0.625); 

\draw[green!50!black, ->-=.7] (-13.75,0.375) -- (-14.5,0); 
\draw[green!50!black, ->-=.7] (-12.5,1) -- (-13.25,0.625); 
\draw[green!50!black, ->-=.7] (-12.75,0.375) -- (-13.5,0); 
\draw[green!50!black, ->-=.7] (-11.5,1) -- (-12.25,0.625); 
\draw[green!50!black, ->-=.7] (-11.75,0.375) -- (-12.5,0); 
\draw[green!50!black, ->-=.7] (-10.5,1) -- (-11.25,0.625); 
\draw[green!50!black, ->-=.7] (-10.75,0.375) -- (-11.5,0); 
\draw[green!50!black, ->-=.7] (-10,0.75) -- (-10.25,0.625); 

\draw[green!50!black, ->-=.7] (-10,0.25) -- (-10.5,0);

\draw[green!50!black] (-15,0.5) -- (-14.7,0.5); 
\draw[green!50!black] (-14.3,0.5) -- (-13.7,0.5); 
\draw[green!50!black] (-13.3,0.5) -- (-12.7,0.5); 
\draw[green!50!black] (-12.3,0.5) -- (-11.7,0.5); 
\draw[green!50!black] (-11.3,0.5) -- (-10.7,0.5); 
\draw[green!50!black] (-10.3,0.5) -- (-10,0.5); 

\draw[green!50!black] (-14.7,0.5) arc (10:-60:0.125); 
\draw[green!50!black] (-13.7,0.5) arc (10:-60:0.125); 
\draw[green!50!black] (-12.7,0.5) arc (10:-60:0.125); 
\draw[green!50!black] (-11.7,0.5) arc (10:-60:0.125); 
\draw[green!50!black] (-10.7,0.5) arc (10:-60:0.125); 
\draw[green!50!black] (-14.3,0.5) arc (190:120:0.125); 
\draw[green!50!black] (-13.3,0.5) arc (190:120:0.125); 
\draw[green!50!black] (-12.3,0.5) arc (190:120:0.125); 
\draw[green!50!black] (-11.3,0.5) arc (190:120:0.125); 
\draw[green!50!black] (-10.3,0.5) arc (190:120:0.125);

        \draw[thick, fill=white] (-15,0) circle(.1); 
        \draw[thick, fill=white] (-14,0) circle(.1); 
        \draw[thick, fill=white] (-13,0) circle(.1); 
        \draw[thick, fill=white] (-12,0) circle(.1); 
        \draw[thick, fill=white] (-11,0) circle(.1); 
        \draw[thick, fill=white] (-10,0) circle(.1); 
        \draw[thick, fill=white] (-15,1) circle(.1); 
        \draw[thick, fill=white] (-14,1) circle(.1); 
        \draw[thick, fill=white] (-13,1) circle(.1); 
        \draw[thick, fill=white] (-12,1) circle(.1); 
        \draw[thick, fill=white] (-11,1) circle(.1); 
        \draw[thick, fill=white] (-10,1) circle(.1); 
        
        \node at (-15,0) {\tiny $A$}; 
        \node at (-14,0) {\tiny $B$}; 
        \node at (-13,0) {\tiny $C$};
        \node at (-12,0) {\tiny $D$}; 
        \node at (-11,0) {\tiny$E$}; 
        \node at (-10,0) {\tiny $A$}; 
        \node at (-15,1) {\tiny $D$}; 
        \node at (-14,1) {\tiny $E$}; 
        \node at (-13,1) {\tiny $A$}; 
        \node at (-12,1) {\tiny $B$}; 
        \node at (-11,1) {\tiny $C$}; 
        \node at (-10,1) {\tiny $D$}; 
    
\end{tikzpicture}
    \caption{$\tau_{M_1} \circ \tau_{M_2} \circ \tau_{M_3} \circ \tau_{M_4} \circ \tau_{M_5}
    (L_1)$}
	\label{fig8}
\end{figure}
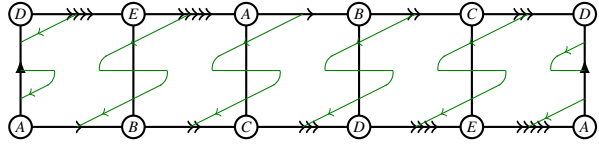 

\begin{figure}[htb!]
\centering
\begin{tikzpicture} [scale=1.5]

        \tikzset{->-/.style={decoration={ markings,
	    mark=at position #1 with {\arrow{>}}},postaction={decorate}}}
        \tikzset{->>-/.style={decoration={ markings,
	    mark=at position #1 with {\arrow{>>}}},postaction={decorate}}}
        \tikzset{->>>-/.style={decoration={ markings,
	    mark=at position #1 with {\arrow{>>>}}},postaction={decorate}}}
       \tikzset{->>>>-/.style={decoration={ markings,
	    mark=at position #1 with {\arrow{>>>>}}},postaction={decorate}}}
        \tikzset{->>>>>-/.style={decoration={ markings,
	    mark=at position #1 with {\arrow{>>>>>}}},postaction={decorate}}}

       \tikzset{-|>-/.style={decoration={ markings,
	    mark=at position #1 with {\arrow{latex}}},postaction={decorate}}}

	\draw[thick, ->-=.55] (-15,0) -- (-14,0);
    \draw[thick, ->>-=.6] (-14,0) -- (-13,0);
    \draw[thick, ->>>-=.65] (-13,0) -- (-12,0);
    \draw[thick, ->>>>-=.7] (-12,0) -- (-11,0);
    \draw[thick, ->>>>>-=.7] (-11,0) -- (-10,0);
        
    \draw[thick, ->>>>-=.65] (-15,1) -- (-14,1);
    \draw[thick, ->>>>>-=.7] (-14,1) -- (-13,1);
    \draw[thick, ->-=.6] (-13,1) -- (-12,1);
    \draw[thick, ->>-=.55] (-12,1) -- (-11,1);
    \draw[thick, ->>>-=.55] (-11,1) -- (-10,1);

	\draw[thick, -|>-=.6] (-15,0) -- (-15,1);
	\draw[thick, -|>-=.6] (-10,0) -- (-10,1);

	\draw[thick] (-14,0) -- (-14,1);
	\draw[thick] (-13,0) -- (-13,1); 
	\draw[thick] (-12,0) -- (-12,1); 
	\draw[thick] (-11,0) -- (-11,1);

\draw[green!50!black, ->-=.7] (-14.5,1) -- (-15,0.75); 
\draw[green!50!black, ->-=.7] (-14.55,0.475) -- (-14.65,0.425); 
\draw[green!50!black, ->-=.7] (-13.5,1) -- (-14.05,0.725); 
\draw[green!50!black, ->-=.7] (-13.55,0.475) -- (-13.65,0.425); 
\draw[green!50!black, ->-=.7] (-14.15,0.175) -- (-14.5,0); 
\draw[green!50!black, ->-=.7] (-12.5,1) -- (-13.05,0.725); 
\draw[green!50!black, ->-=.7] (-12.55,0.475) -- (-12.65,0.425); 
\draw[green!50!black, ->-=.7] (-13.15,0.175) -- (-13.5,0); 
\draw[green!50!black, ->-=.7] (-11.5,1) -- (-12.05,0.725); 
\draw[green!50!black, ->-=.7] (-11.55,0.475) -- (-11.65,0.425); 
\draw[green!50!black, ->-=.7] (-12.15,0.175) -- (-12.5,0); 
\draw[green!50!black, ->-=.7] (-10.5,1) -- (-11.05,0.725); 
\draw[green!50!black, ->-=.7] (-10.55,0.475) -- (-10.65,0.425); 
\draw[green!50!black, ->-=.7] (-11.15,0.175) -- (-11.5,0); 
\draw[green!50!black, ->-=.7] (-10,0.75) -- (-10.05,0.725); 
\draw[green!50!black, ->-=.7] (-10.15,0.175) -- (-10.5,0);

\draw[green!50!black] (-15,0.6) -- (-14.5,0.6); 
\draw[green!50!black] (-14.1,0.6) -- (-13.5,0.6); 
\draw[green!50!black] (-13.1,0.6) -- (-12.5,0.6); 
\draw[green!50!black] (-12.1,0.6) -- (-11.5,0.6); 
\draw[green!50!black] (-11.1,0.6) -- (-10.5,0.6); 
\draw[green!50!black] (-10.1,0.6) -- (-10,0.6); 

\draw[green!50!black] (-14.7,0.3) -- (-14.1,0.3); 
\draw[green!50!black] (-13.7,0.3) -- (-13.1,0.3); 
\draw[green!50!black] (-12.7,0.3) -- (-12.1,0.3); 
\draw[green!50!black] (-11.7,0.3) -- (-11.1,0.3); 
\draw[green!50!black] (-10.7,0.3) -- (-10.1,0.3); 

\draw[green!50!black] (-14.7,0.3) arc (190:120:0.125); 
\draw[green!50!black] (-13.7,0.3) arc (190:120:0.125); 
\draw[green!50!black] (-12.7,0.3) arc (190:120:0.125); 
\draw[green!50!black] (-11.7,0.3) arc (190:120:0.125); 
\draw[green!50!black] (-10.7,0.3) arc (190:120:0.125); 

\draw[green!50!black] (-14.1,0.3) arc (10:-60:0.125); 
\draw[green!50!black] (-13.1,0.3) arc (10:-60:0.125); 
\draw[green!50!black] (-12.1,0.3) arc (10:-60:0.125); 
\draw[green!50!black] (-11.1,0.3) arc (10:-60:0.125); 
\draw[green!50!black] (-10.1,0.3) arc (10:-60:0.125); 

\draw[green!50!black] (-14.5,0.6) arc (10:-60:0.125); 
\draw[green!50!black] (-13.5,0.6) arc (10:-60:0.125); 
\draw[green!50!black] (-12.5,0.6) arc (10:-60:0.125); 
\draw[green!50!black] (-11.5,0.6) arc (10:-60:0.125); 
\draw[green!50!black] (-10.5,0.6) arc (10:-60:0.125); 

\draw[green!50!black] (-14.1,0.6) arc (190:120:0.125); 
\draw[green!50!black] (-13.1,0.6) arc (190:120:0.125); 
\draw[green!50!black] (-12.1,0.6) arc (190:120:0.125); 
\draw[green!50!black] (-11.1,0.6) arc (190:120:0.125); 
\draw[green!50!black] (-10.1,0.6) arc (190:120:0.125);

        \draw[thick, fill=white] (-15,0) circle(.1); 
        \draw[thick, fill=white] (-14,0) circle(.1); 
        \draw[thick, fill=white] (-13,0) circle(.1); 
        \draw[thick, fill=white] (-12,0) circle(.1); 
        \draw[thick, fill=white] (-11,0) circle(.1); 
        \draw[thick, fill=white] (-10,0) circle(.1); 
        \draw[thick, fill=white] (-15,1) circle(.1); 
        \draw[thick, fill=white] (-14,1) circle(.1); 
        \draw[thick, fill=white] (-13,1) circle(.1); 
        \draw[thick, fill=white] (-12,1) circle(.1); 
        \draw[thick, fill=white] (-11,1) circle(.1); 
        \draw[thick, fill=white] (-10,1) circle(.1); 
        
        \node at (-15,0) {\tiny $A$}; 
        \node at (-14,0) {\tiny $B$}; 
        \node at (-13,0) {\tiny $C$};
        \node at (-12,0) {\tiny $D$}; 
        \node at (-11,0) {\tiny$E$}; 
        \node at (-10,0) {\tiny $A$}; 
        \node at (-15,1) {\tiny $D$}; 
        \node at (-14,1) {\tiny $E$}; 
        \node at (-13,1) {\tiny $A$}; 
        \node at (-12,1) {\tiny $B$}; 
        \node at (-11,1) {\tiny $C$}; 
        \node at (-10,1) {\tiny $D$}; 
    
\end{tikzpicture}
    \caption{$Cone(Cone((M_1 \oplus M_2 \oplus M_3 \oplus M_4 \oplus M_5) \to L_1) \to L_1) \simeq
    L_2$}
	\label{fig9}
\end{figure}
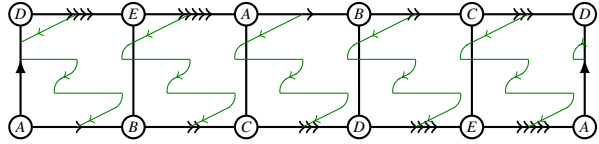 

When $A \neq \emptyset$, if the puncture between $M_i$ and $M_{i+1}$ is closed, then they become
isotopic hence give equivalent objects. (Of course, one has to isotope them with a finger move so as
to make both of them exact Lagrangians). The same argument as above, with the understanding
that some of the $M_i$ represent equivalent objects, shows the exact triangle
(\ref{eq:exactness}) remains valid. Hence, $L_1$ and $L_2$ again split generates
$\mathcal{F}(\Sigma)$.  \ed

\section{Application to vector bundles over cycles over projective lines}\label{cycles-sec}

\subsection{Simple vector bundles on cycles of projective lines}

In this subsection we work over an algebraically closed field $k$ of characteristic $\neq 2$.
Let $C=\cup_{j=0}^{n-1} C_j$ be a cycle of $n$ projective lines (also known as the standard $n$-gon).
We identify each $C_j$ with the standard copy $\P^1$ in such a way that the point $\infty\in C_j$ is glued
to the point $0\in C_{j+1}$ (we identify indices with $\Z/n$).

Recall that, up to isomorphism, all simple vector bundles on $C$ are obtained by the following construction (see \cite{BDG}),
which has as an input an integer valued matrix $\bm=(m^j_i)_{i=1,\ldots,r;j=0,\ldots,n-1}$ and a nonzero constant $\la\in k^*$.
The corresponding vector bundle $V=V^\la(\bm)$ is defined by setting
$$V|_{C_j}=V_j=\OO_{\P^1}(m^j_1)\oplus\ldots\oplus \OO_{\P^1}(m^j_r)$$
and by making the following identifications $V_j|_{\infty}\simeq V_{j+1}|_0$: for all $j$ except for $j=n-1$
we use the standard trivializations of the corresponding bundle $\OO(m)$ at $0$ and at $\infty$ 
  (given by $x_0^m$ and $x_1^m$), 
  while for $j=n-1$, we
use
$$\la\cdot C: V_{n-1}|_{\infty}\to V_0|_0,$$
where $C$ is the transitive permutation matrix $e_i\mapsto e_{i-1}$ (where the indices are in $\Z/r\Z$. 
The obtained vector bundle of rank $r$ is simple
if and only if a certain condition on $\bm$ is satisfied. Namely, let us unroll the matrix $\bm$ into an $rn$-periodic sequence
by setting 
$$d_{qn+j}=m^j_{-q}, \ \ j, q\in \Z, 0\le j<n.$$
The conditions are: (1) for every $i,i',j$ one has $|m^j_i-m^j_{i'}|\le 1$;
(2) for every $q$, not divisible by $r$, the $rn$-periodic sequence $(d_{qn+j}-d_j)$ is not identically $0$ and the occurrences
of $1$ and $-1$ in it alternate.

Recall that one of the results of \cite{pol09} is an explicit computation of 
the trigonometric solution of the AYBE associated with a pair $(V,\OO_p)$, where $V=V^\la(\bm)$ is a simple
bundle on $C$ and $p$ is a smooth point. The answer is given by the trigonometric solution corresponding
to a certain associative Belavin-Drinfeld structure $\ABD(V,p)$, which we will describe now.
Without loss of generality we can assume that $p\in C_0=C_n$. Let us define the complete order $\prec$ on the set of indices
$\Z/r\Z=\{0,1,\ldots,r-1\}$ as follows: $i\prec i'$ if the sequence $(d_{j-in}-d_{j-i'n})_{j=0,1,\ldots}$ is 
                                                 nonzero and the first nonzero term in it is
                                                 negative
(the fact that it is a complete order follows from the condition (2) above).
We define the transitive permutations $C_1$ and $C_2$ on $\Z/r\Z$ by letting $C_1$ send each non-maximal element with
respect to the above complete order to the next element, and by $C_2(i)=i-1$. Finally, we define a subset 
$A\sub \Z/r\Z$ to be the set of $i$ such that $i-1\prec C_1(i)-1$ and $m^j_i=m^j_{i'}$ for $0<j<n$.
By \cite[Thm.\ 5.3]{pol09}, in fact $C_2$ is a power of $C_1$, and the solution of the AYBE associated with a natural family of deformations of $V$ and $p$ is  the solution \eqref{trig-AYBE-eq} associated with
$$\ABD(V,p):=(\Z/r,C_1,C_2,A).$$ 

Now the arguments of Section \ref{alg-form-AYBE-sec} imply that the formal solution of the general AYBE associated with
the pair $(V,p)$ is equivalent to \eqref{trig-AYBE-eq}, viewed as a formal solution.
By Theorem A, this implies that the $A_\infty$-subcategory, split generated by $V$ and $\OO_p$, depends only on 
$\ABD(V,p)$. Here to apply Theorem A (with $R=k$) we need to equip the $A_\infty$-algebra of endomorphisms
  of $V\oplus\OO_p$ with a cyclic structure with respect to a natural pairing coming from the Serre duality. 
  The existence of such a cyclic structure can be proved similarly to \cite[Sec.\ 4.8]{pol18} (using the assumption that 
  characteristic is not equal to $2$). Namely, first,
  using a $1$-spherical twist we can replace $V\oplus \OO_p$ with a vector bundle, and then, use
  Proposition 4.8.2 and Lemma 4.8.4 of \cite{pol18}. In the characteristic zero case one can instead use the criterion 
  of Kontsevich-Soibelman \cite[Thm.\ 10.2.2]{KS} (see \cite[Rem.\ 4.8.3]{pol18}).

%Discuss generation for pairs $(\OO,V)$ (resp., $(\OO_p,V)$), where $V$ is a vector bundle on $C$.

\begin{defi} We say that a vector bundle $W$ on $\P^1$ is of {\it positive (resp., nonnegative) type}
if $W\simeq\bigoplus_{j=1}^r \OO_{\P^1}(a_i)$ with all $a_i>0$ (resp., $a_i\ge 0$).
Now let $V$ be a vector bundle on $C$. We say that $V$ is of {\it positive (resp., nonnegative) type} if
each restriction $V|_{C_i}$ is of positive (resp., nonnegative) type.
In the case $n=1$ we require this property for $f^*V$, where $f:\P^1\to C$ is the normalization map.
\end{defi}
  
  Recall that a collection of objects $(O_i)$ {\it split generates} a triangulated category $\TT$ if the minimal
  triangulated subcategory $\TT'\sub \TT$, closed under direct summands and containing all $O_i$, is the entire $\TT$.

\begin{prop}\label{gen-prop} 
Let $V$ be a simple vector bundle on $C$ of positive type. Then the pair $(\OO_C,V)$ split generates
$\Perf(C)$.
\end{prop}

\Pf . Let us pick smooth points $p_1,\ldots,p_n$, one on each component of $C$.
Note that $V(-p_1-\ldots-p_n)$ is of nonnegative type. 

Assume first that the rank of $V$ is $>1$.
Hence, by Lemma \ref{subbun-lem}(i) below, there exists an injection of $\OO_C$ into $V(-p_1-\ldots-p_n)$.
Let us consider the composed injective morphism
$$f:\OO_C\to \OO_C(p_1+\ldots+p_n)\to V,$$
where the first arrow is given by the canonical section of $\OO_C(p_1+\ldots+p_n)$ vanishing at 
the divisor $p_1+\ldots+p_n$. Then the coherent sheaf $\coker(f)$ has nonzero torsion at each of the points
$p_1,\ldots,p_n$. 
Thus, we have an exact sequence
$$0\to \bigoplus_{i=1}^n \TT_i\to \coker(f)\to \FF\to 0$$
where $\TT_i$ is a nonzero sheaf supported at $p_i$ and $\FF$ is locally free near each $p_1,\ldots,p_n$.
Such a sequence necessarily splits, so each $\TT_i$ is a direct summand of $\coker(f)$.
This shows that the subcategory, split generated by $\OO_C$ and $V$ contains $\TT_1,\ldots,\TT_n$.
Furthermore, each $\TT_i$ has a direct summand of the form $\OO_{m_ip_i}$ with some $m_i\ge 1$.
It remains to note that the objects $(\OO_C,\OO_{m_1p_1},\ldots,\OO_{m_np_n})$ split generate $\Perf(C)$.
Indeed, this can be checked similarly to \cite[Lem.\ 3.3.1]{LP}: starting from $\OO_C$ and using the exact sequences
of the form 
$$0\to L(-m_ip_i)\to L\to \OO_{m_ip_i}\to 0,$$
we derive that all the negative powers of the ample line bundle $\OO_C(\sum m_ip_i)$ belong to the subcategory split generated
by our objects. The fact that all negative powers of an ample line bundle generate $\Perf(C)$ is proved in
  \cite[Thm.\ 4]{Orlov-gen}.

In the case when $V$ is a line bundle, of positive degree on each component, by Lemma \ref{subbun-lem}(i), 
we can find a global section $s:\OO_C\to V$ which
does not vanish at the nodes. Its restriction to every component of $C$ vanishes at some smooth point $p_i$.
Then $\coker(s)$ will again have a nonzero torsion part at each $p_i$, and the above proof goes through. 
\ed

\begin{lem}\label{subbun-lem} 
(i) Let $W$ be a simple vector bundle on $C$ of nonnegative type. Assume in addition that either $W$ has rank $>1$,
or has positive degree. Then there exists an injective morphism
$\OO_C\to W$, which is an embedding as a subbundle near the nodes.

  \noindent (ii) Let $V$ be a simple vector bundle of positive type, $p\in C$ a smooth point. 
Let us denote by $E(V,p)$ the universal extension 
$$0\to \Ext^1(V,\OO_C(p))^*\ot \OO_C(p)\to E(V,p)\to V\to 0.$$
Then $E(V,p)$ is the result of applying to $V$ the inverse twist with respect to $\OO_C(p)$.
In particular, $E(V,p)$ is still a simple vector bundle.  
  \end{lem}

\Pf . (i) We use the fact that $W$ has the form $W=V^\la(\bm)$, where all $m^j_i\ge 0$. Note that the condition that $W$
is simple and has rank $>1$ implies that $m^j_i>0$ for at least one pair $(i,j)$. To define a global section of
$W$ we need to choose a global section $s^j_i\in H^0(C_j,\OO(m^j_i)$ for each $(i,j)$ in such a way that they are compatible
with the gluing over $C_j\cap C_{j+1}$. We claim that we can make these choices in such a way that each $s^j_i$ is nonzero
at $0$ and $\infty$. Indeed, in the case when $m^j_i>0$ we can arrange $s^j_i$ to have arbitrary values at $0$ and $\infty$,
while in the case $m^j_i=0$, one of these values determine the other. Now looking at the way the gluing is defined for 
$V^\la(\bm)$ we see that the existence of at least one positive $m^j_i$ guarantees the existence of
a global section which is nonzero at all the nodes.

  \noindent
  (ii) This follows from the vanishing $\Hom(V,\OO_C(p))^*\simeq H^1(V(-p))=0$ 
  that holds since $V(-p)$ is of nonnegative type. Note that by Serre duality,
  any line bundle on $C$ is $1$-spherical. 
  \ed

Now we are going to consider the associative Belavin-Drinfeld structure  $\ABD(E(V,p),p)$
 associated to $E(V,p)$ and $p$.

\begin{thm}\label{ABD-equiv-thm} 
Let $V$ and $V'$ be simple vector bundles on $C$ of positive type.
Assume that for some smooth points $p,p'\in C$ one has an isomorphism
$$\ABD(E(V,p),p)\simeq \ABD(E(V',p'),p')$$ 
of associative Belavin-Drinfeld structures. Then there exists a Fourier-Mukai
autoequivalence $\Phi$ of $\Perf(C)$ given by a kernel in $D^b(C\times C)$, such that
$\Phi(\OO_C)\simeq \OO_C$ and $\Phi(V)\simeq V'$.
\end{thm}

\Pf . By Lemma \ref{subbun-lem}(ii), the inverse twist with respect to $\OO_C(p)$ sends the pair $(\OO_C,V)$ to the pair $(O_p,E(V,p)[1])$.
Similarly, the twist with respect to $\OO_C(p')$ sends $(\OO_C,V')$ to $(\OO_{p'},E(V,p)[1])$. By Theorem A,
the isomorphism of the corresponding associative Belavin-Drinfeld structures implies that the subcategories,
split generated by $(\OO_C,V)$ and $(\OO_C,V')$ are related by an equivalence $\Phi$ in such a way that 
  $\Phi(\OO_C)\simeq \OO_C$
and $\Phi(V)\simeq V'$. By Proposition \ref{gen-prop}, $\Phi$ is actually an autoequivalence of $\Perf(C)$, or more precisely,
of its $A_\infty$-enhancement. Such an autoequivalence is always given by a kernel on $C\times C$ which could be a complex
of quasicoherent sheaves (see \cite{Toen}). 
The fact that it belongs to the bounded derived category of coherent sheaves follows from \cite[Lem.\ 3.5.1]{LP}.
\ed

\subsection{Proof of Theorem C}\label{proof-C-sec}

It is enough to consider the case when $V$ is of positive type. Indeed, starting from an arbitrary bundle we can apply
twists at smooth points to replace $V$ with $V(N(p_1+\ldots+p_n))$, which is of positive type for large $N$.

By Lemma \ref{subbun-lem}(ii), the inverse twist with respect to $\OO_C(p)$ transforms the pair $(\OO_C,V)$ to
$(\OO_p,E(V,p)[1])$. As was shown in \cite{pol09}, the solution of AYBE, associated with the pair $(E(V,p),\OO_p)$,
is a trigonometric solution \eqref{trig-AYBE-eq}, corresponding to an associative Belavin-Drinfeld structure $(S,C_1,C_2,A)$
in which $C_2=C_1^k$ for some $k$. Hence, by Theorem A, Theorem B and
Proposition \ref{g1-Fuk-gen-prop}, the subcategory in $D^b(C)$ split generated by the pair $(\OO_C,V)$ is equivalent to the 
                                                 Fukaya category of some square-tiled
surface of genus $1$, in such a way that $\OO_C$ and $V$ correspond to the Lagrangians $L_1$ and $L_2$. Note that in establishing this
equivalence we apply Theorem A, so we pass to formal solutions of the AYBE, as explained in Section \ref{alg-form-AYBE-sec}.

Now we recall that the Dehn twists with respect to the graded Lagrangians $L_1,M_1,\ldots,M_n$ generate the pure mapping class
                                                 group (see the proof of Proposition \ref{g1-Fuk-gen-prop}).
                                                 Hence, there exists a composition
of these Dehn twists and their inverses that takes $L_1$ into $L_2$. 
  Under the above equivalence, this corresponds to a composition $\Phi$
of $1$-spherical twists and their inverses that takes $\OO_C$ into $V$.
\ed

\end{document}